\documentclass[12pt]{amsart}

\makeatletter
\renewcommand{\subsection}{\@startsection
{subsection}%
{2}%
{0mm}%
{-\baselineskip}%
{-0.5\baselineskip}%
{\normalfont\normalsize\bfseries }}%
\makeatother

\usepackage[all]{xy}
\usepackage{amssymb,amscd}
\usepackage{url}

 \usepackage{hyperref}

\newtheorem{thm}{Theorem}
\newtheorem{lemma}[thm]{Lemma}
\newtheorem{cor}[thm]{Corollary}
\newtheorem{prop}[thm]{Proposition}

\theoremstyle{definition}

\newtheorem{rem}[thm]{Remark}

\theoremstyle{remark}

\newcommand{\ce}{\mathcal{E}}
\newcommand{\cf}{\mathcal{F}}

\newcommand{\cl}{\mathcal{L}}
\newcommand{\cR}{\mathcal{R}}
\newcommand{\cs}{\mathcal{S}}
\newcommand{\ct}{\mathcal{T}}

\newcommand{\al}{\alpha}
\newcommand{\lb}{\lambda}
\newcommand{\Lb}{\Lambda}

\newcommand{\de}{\delta}
\newcommand{\De}{\Delta}

\newcommand{\w}{\omega}

\newcommand{\W}{\Omega}

\newcommand{\g}{\gamma}
\newcommand{\ve}{\varepsilon}
\newcommand{\G}{\Gamma}

\newcommand{\s}{\sigma}
\newcommand{\na}{\nabla}

 \newcommand{\ify}{\infty}

 \newcommand{\id}{\mathrm{id}}

 \newcommand{\Ab}{\mathbb{A}}
 \newcommand{\Cb}{\mathbb{C}}
 \newcommand{\Rb}{\mathbb{R}}
 
 \newcommand{\Zb}{\mathbb{Z}}

 \DeclareMathOperator{\ad}{ad}

 \DeclareMathOperator{\ord}{order}
 \DeclareMathOperator{\Ind}{Ind}

 \DeclareMathOperator{\End}{End}
 \DeclareMathOperator{\Diff}{Diff}
 \DeclareMathOperator{\Ch}{Ch}
 
 \DeclareMathOperator{\Tr}{Tr}
 
 \DeclareMathOperator{\Hom}{Hom}

\begin{document}
\title{ Bivariant Chern Character and the Longitudinal Index Theory.}
\author{Alexander Gorokhovsky}

\address{Department of Mathematics\\
University of Colorado\\UCB 395\\ Boulder, Colorado 80309-0395\\
USA}

\email {Alexander.Gorokhovsky@colorado.edu}

\thanks{The author thanks A.~Connes, X.~Dai, J.~Lott, R.~Nest, V.~Nistor, H.~Moscovici, P.~Piazza,
R.~Ponge and B.~Tsygan for many helpful discussions. }

\begin{abstract} In this paper we consider a family of Dirac-type
operators on fibration $P \to B$ equivariant with respect to an
action of an \'etale groupoid.  Such a family defines an element in
the bivariant $K$ theory. We compute the action of the bivariant
Chern character of this element on the image of Connes' map $\Phi$
in the cyclic cohomology. A particular case of this result is
Connes' index theorem for \'etale groupoids \cite{cn94} in the case
of fibrations.
\end{abstract}

 \maketitle

\section{Introduction}\label{intro}

   In \cite{ni91, ni93a}, answering a question posed by
A.~Connes in \cite{cn83}, V.~Nistor defined a bivariant Chern
character for a $p$-summable quasihomomorphism and established its
fundamental properties. This theory was further developed in
\cite{cq95a, cuntz97}. Bivariant Chern character encodes extensive
information related to index theory. In the present paper we compute
the action of the bivariant Chern character on cyclic cohomology in
geometric situations arising from equivariant families of elliptic
operators.

We consider first the following geometric situation. Let $\pi:P \to
B$ be a fibration. We assume that we are given a vertical Riemannian
metric on this fibration, i.e. that there is a Riemannian metric on
each of the fibers $P_b =\pi^{-1}(b)$, which varies smoothly with
$b$. We assume that each $P_b$ is a complete Riemannian manifold.
 Let $D$ be a family of fiberwise
Dirac type operators on this fibration acting on the section of the
bundle $\ce$. In this paper we consider only the even-dimensional
situation, and hence the bundle $\ce$ is always $\Zb_2$-graded. We
assume that an \'etale groupoid $G$ with the unit space $B$ acts on
this fibration, and that the operator $D$ is invariant under this
action. This means that  every $\g \in G$ defines a diffeomorphism
$P_{r(\g)} \to P_{s(\g)}$, $p \mapsto p\g$.One also requires that
these diffeomorphisms are compatible with the groupoid structure,
i.e. that $(p \g_1)\g_2 =p(\g_1\g_2)$. Notice that this implies that
the vertical metric is invariant under the action of $G$. We do not
assume existence of a metric on $P$, invariant under the action of
$G$.

There are two constructions in $K$-theory and cyclic cohomology we
need to use to describe the problem. The first is Nistor's bivariant
Chern character. It appears in our setting as follows. The operator
$D$ described above  defines an element of the equivariant
$KK$-theory $KK^{G}\left( C_0(P), C_0(B)\right)$, \cite{lg99a,
lg99b}. Using the canonical map
\begin{equation}
j^{G}: KK^{G}\left( C_0(P), C_0(B)\right) \to KK\left( C_0(P)
\rtimes G, C_0(B)\rtimes G\right)
\end{equation}
we obtain a class in $ KK\left( C_0(P) \rtimes G, C_0(B)\rtimes
G\right)$ defined by $D$.  This class can be represented by an
explicit quasihomomorphism $\psi_D$ in the sense of J.~Cuntz
\cite{cuntz81, cuntz82}. Moreover, one can show that $\psi_D$
actually defines a \emph{p-summable} quasihomomorphism
 of smooth algebraic cross-products $C_0^{\ify}(P) \rtimes
G$ and $C_0^{\ify}(B) \rtimes G$ in the sense of V.~Nistor
\cite{ni91, ni93a}. One can then use  techniques developed in
\cite{ni91, ni93a}  to define the bivariant Chern character $
\Ch(D)$ in bivariant cyclic homology which is a morphism of
complexes $CC_*(C_0^{\ify}(P) \rtimes G) \to CC_*(C_0^{\ify}(B)
\rtimes G)$.

The other tool which we need  is the explicit construction, due to
A.~Connes \cite{cn86, cn94}, of the classes in the cyclic cohomology
of cross product algebras. Let $M$ be a  manifold on which an
\'etale groupoid $G$ acts, and let $M_{G}$ be the corresponding
homotopy quotient.  Then Connes constructs an explicit chain map of
complexes inducing an injective map in cohomology $\Phi:
H^*_{\tau}\left(M_{G}\right) \to HC^*\left(C_0^{\ify}(M)\rtimes G
\right)$. Here $\tau$ denotes a twisting by the orientation bundle
of $M$. Introduce now the following notations. The map $\pi$ induces
a map $P_G \to B_G=BG$ which we also denote by $\pi$.  Since the
fibers of $P_G \to B_G$ are $ spin^c$ and hence oriented we have a
pull-back map $\pi^*: H^*_{\tau}\left(B_{G}\right) \to
H^*_{\tau}\left(P_{G}\right)$. Let $\widehat{A}_G(TP/B) \in
H^*(P_G)$ be the equivariant $\widehat{A}$-genus of the vertical
tangent bundle, i.e. the $\widehat{A}$-genus   of the bundle on
$P_G$ induced by the vertical tangent bundle  $TP/B$ on $P$. Note
that we use the conventions from \cite{bgv} in the definitions of
characteristic classes. Let $Ch_G(\ce/\cs)$ be the equivariant
twisting Chern character of the bundle $\ce$. If the fibers of $P
\to B$ have $spin$ structure, and $\ce =\cs \otimes V$ where $\cs$
is the vertical spin bundle, then $\Ch_G(\ce/\cs)=\Ch_G(V)$ is just
the equivariant Chern character of $V$. In other words it is the
Chern character of the bundle on $P_G$ induced by the  bundle $V$ on
$P$. Set now for $c \in H^*_{\tau}\left(B_{G}\right)$
$\widetilde{\pi}^*(c) =(2\pi i)^{-\frac{\dim P-\dim B}{2}}
\widehat{A}_G(TP/B) \Ch_G(\ce/\cs) \pi^*(c)$

 Our main result is then the following:
\begin{thm} \label{mainthm}
The diagram
\begin{equation} \label{comm}
 \xymatrix{& H^*_{\tau}(P_{G}) \ar[rr]^{\Phi \ \ \ \ \ \ } & &HC^* ( C_0^{\ify}(P)\rtimes G )\\
 &&&\\
 &\ar[uu]^{\widetilde{\pi}^*} H^*_{\tau}(B_{G}) \ar[rr]^{\Phi \ \ \ \ \ \
}& & HC^* ( C_0^{\ify}(B)\rtimes G )\ar[uu]_{\Ch(D)^t}
         }
\end{equation}
commutes.
\end{thm}

As an illustration   consider the case when  the action of $G$ on
$P$ is free, proper and cocompact. In this case there is a Morita
equivalence between the algebras $C_0^{\ify}(P)\rtimes G$ and
$C_0^{\ify}(P/G)$. We have a canonical class $1 \in
K_0\left(C_0^{\ify}(P/G)\right)$, and hence the corresponding class
in $K_0\left(C_0^{\ify}(P)\rtimes G\right)$. This class can be
represented   by an idempotent $e\in C_0^{\ify}(P)\rtimes G $. It
can be described explicitly as follows. Let $\phi \in C_0^{\ify}(P)$
be such that for every $p \in P$ $\sum_{r(\g)=\pi(p)}
\phi(p\g)^2=1$. We then define an element in $C_0^{\ify}(P) \rtimes
G$ by $e(p,\g)=\phi(p) \phi(p \g)$. Since the action of $G$ on $P$
is proper one can define and index $\Ind D$ of $D$ in
$K_0(C_0^{\ify}(G)\otimes \cR)$ where $\cR$ is the algebra of
matrices with rapidly decaying entries \cite{cm90}. One can
represent index of $D$ as the product in $KK$-theory: $\Ind D
=\psi_D [e]$. This, together with commutativity of \eqref{comm}
allows to compute the pairing of $\Ind D$ with $H^*(B_{G})$ in
topological terms:
\begin{multline}
\left\langle
 \Phi(c), \Ch\left(\Ind D \right) \right \rangle=\\
(2\pi i)^{-\frac{\dim P-\dim B}{2}}\left \langle
\Phi(\widehat{A}_G(TP/B) \Ch_G(\ce/\cs) \pi^*(c)), \Ch e \right
\rangle
\end{multline}
Now due to the conditions imposed on the action in this case $P/G$
is a smooth manifold, and every cocycle $C \in C^*(G, \W_*(P))$
defines a class $[C]\in H^*_{\tau}(P/G)$. It is easy to see that
\begin{equation}
\left \langle\Phi(C), \Ch e \right \rangle =\int \limits_{P/G}[C].
\end{equation}
Hence the pairing in the right hand side equals $\int \limits_{P/G}
\widehat{A}(\cf) \Ch(\ce/\cs) [\pi^*(c)]$, where $\cf$ is the
foliation on $P/G$ induced by the fibers of the fibration $\pi$.
Thus in this case we recover Connes' index theorem for \'etale
groupoids \cite{cn94}, for the case when $P\to B$ is a fibration.
Notice that we have a constant different from the one in
\cite{cn94}, which is due to the different conventions in defining
characteristic classes. Among important particular cases of this
theorem let us mention Connes'-Moscovici higher $\G$-index theorem,
obtained when the groupoid $G$ is a discrete group. There are
several other proofs of this theorem using the bivariant Chern
character idea, compare \cite{wu97, perd04}.

Our proof or the Theorem \ref{mainthm} is as follows. Let $\W_*(B)$
denote the complex of smooth currents, i.e. differential forms with
values in the orientation bundle, on the manifold $B$. Then
cohomology of the complex $C^*(G, \W_*(B))$ equals
$H^*_{\tau}(B_G)$. We show that the diagram of complexes
corresponding to \eqref{comm} commutes up to the chain homotopy of
complexes. To this end we use a simplicial version $\Ab$
\cite{dup76, dup78} of the Bismut superconnection \cite{bis85} to
construct a map of complexes $\Phi_{\Ab}: C^*(G, \W_*(B)) \to
CC^*(C_0^{\ify}(P)\rtimes G)$ to obtain the following diagram:
\begin{equation}\label{comm'}
 \xymatrix{& C^*(G, \W_*(P)) \ar[rr]^{\Phi \ \ \ \ \ \ } & &CC^* ( C_0^{\ify}(P)\rtimes G )\\
 &&&\\
 &\ar[uu]^{\widetilde{\pi}^*} C^*(G, \W_*(B))\ar@{-->}[uurr]^{\Phi_{\Ab}} \ar[rr]^{\Phi \ \ \ \ \ \
}& & CC^* ( C_0^{\ify}(B)\rtimes G )\ar[uu]_{\Ch(D)^t}
         }
\end{equation}

We then show that each of the triangles in this diagram is
commutative up to homotopy. The map $\Phi_{\Ab}$ plays the role of
McKean-Singer formula in our context. To show the commutativity of
the upper triangle we replace the superconnection $\Ab$ by a
rescaled superconnection $\Ab_s$ and compute the limit when $s \to
0$. To show commutativity of the lower triangle the classical method
\cite{bis85} of computing the limit when $s \to \infty$ runs into
serious difficulties, as explained in \cite{ni95}. We avoid this
difficulty by using a modification of the method from \cite{gor03,
gor05}.

We note that superconnection proofs of the Connes-Moscovici higher
$\G$-index theorem were obtained in \cite{lott92, wu97}.

The paper is organized as follows. In the Section \ref{construction}
we construct the map $\Phi_{\Ab}$. Then, after some preparations in
the Section \ref{fpdo}, we prove commutativity of the lower triangle
in the section \ref{biv}. In the Section \ref{bismut} we use Bismut
superconnection to show commutativity of the upper triangle.

\section{Construction of the map $\Phi_{\Ab}$.}\label{construction}

In this section we construct the map $\Phi_{\Ab}$, as described in
the introduction. We start by considering in \ref{trivial} the case
when our groupoid is just an ordinary manifold. Then, after
reviewing in \ref{groupoid} some definitions and results about
actions of groupoids on algebras, we proceed to give the general
construction in \ref{general}.

\subsection{}{\label{trivial}
Consider first the following situation. Let $P \to B$ be a
submersion. We use the notation from \cite{gj93} for the cyclic
complexes. For an algebra $A$ set $C^k\left(A\right)= \left(A
\otimes \bar{A}^{\otimes k}\right)'$, where $\bar{A} =A/\Cb 1$.  Let
$u$ be a formal variable of degree $2$. We denote by
$CC^*\left(A\right)$ the periodic complex of $A$: the complex
$\left(C^*\left(A\right)[u,u^{-1}], b+uB\right)$. For the nonunital
algebras  we consider the reduced cyclic complex of the
unitalization. Another complex we consider is the complex of smooth
currents. We use notation $\W_k \subset \left(\W^k\right)'$ for
smooth currents of degree $k$. We also adjoin $u$ to this complex,
so that $\left(\W_*\right)^k =\oplus_{i} u^i \W_{k-2i}$, and the
differential of degree one is given by $u
\partial$ where $\partial =d^t$ is the transpose of de Rham
differential.
 In this
situation we consider the complex $ \Hom \left(\W_*\left(B\right),
CC^*\left(C_0^{\ify}\left(P\right)\right)\right)$. Here we consider
complex of homomorphisms of $\mathbb{C}[u, u^{-1}]$ modules. Let $D$
be a family of Dirac-type operators on this submersion, acting on
the sections of a bundle $\ce$. On the base $B$ we have an
infinite-dimensional bundle $\pi_*(\ce)$ whose fiber over a point $b
\in B$ is $\G\left(\ce, P_b \right)$. Consider now a superconnection
$\Ab$ adapted to the operator $D$, compare \cite{bgv}. We assume
that the superconnection has the form $\Ab=D+\Ab_{[1]}+\ldots$ where
 $\Ab_{[1]}$ is a connection on the
bundle $\pi_*\left(\ce\right)$ and $\Ab_{[i]}$ for $i> 1$ is a
proper fiberwise pseudodifferential operator.  We assume that the
connection in a local trivialization chart has a form $d + \w$ where
$\w$ is a $1$-form on $B$ with values in the proper fiberwise
pseudodifferential operators. With such a superconnection we will
associate a cochain $\Theta_{\Ab}$ of degree $0$ in the complex $
\Hom \left(\W_*\left(B\right),
CC^*\left(C_0^{\ify}\left(P\right)\right)\right)$.

As a first step we  have the following proposition.

\begin{prop}\label{cocycle}
The following expression defines a cocycle of degree $0$ in the
complex
 $\Hom \left(\W_*\left(B\right),
C^*\left(C_0^{\ify}\left(P\right)\right)[u, u^{-1}]], b+uB\right)$

\begin{equation}
\theta_{\Ab}: c \mapsto \sum \limits_{l \ge  -\frac{\deg c}{2} }
u^{-l}\theta_l\left(c\right)
\end{equation}
where $\theta_l \in C^k\left(C_0^{\ify}\left(P\right)\right)$, $k=
\deg c +2l$, is given by the formula:
\begin{multline}\label{defj}
\theta_l\left(c\right) \left(a_0, a_1, \ldots , a_k\right)\\= \left
\langle c, \int \limits_{\De^k} \Tr_s a_0 e^{-t_0\Ab^2} [\Ab, a_1]
\ldots [\Ab, a_k] e^{-t_k\Ab^2} dt_1\ldots dt_k\right \rangle
\end{multline}
Here $\De^k = \{\left(t_0, t_1, \ldots t_k\right) \ | \ \sum t_i=1,
t_i \ge 0\}$.
\end{prop}
\begin{proof}
It is easy to see that the expression under the trace is a always a
smoothing operator, and hence the trace is well defined. It is also
easy to see that the expression in the formula \eqref{defj} is
nonzero only if $\deg c -k$ is even. The cocycle condition is
verified by a direct computation compare \cite{gs89}.
\end{proof}

We can now follow the method of Connes and Moscovici \cite{cm93} and
replace $\theta_{\Ab}$ by a  finite  cochain $\Theta_{\Ab} \in \Hom
\left(\W_*\left(B\right),
CC^*\left(C_0^{\ify}\left(P\right)\right)\right)$.

The construction is as follows. Let $\Ab_s$ denote the rescaled
superconnection $\Ab_s=sD+\Ab_{[0]}+s^{-1}\Ab_{[1]}+\ldots$.

 Introduce the cochain $\tau_s \in \Hom^{-1} \left(\W_*\left(B\right),
C^*\left(C_0^{\ify}\left(P\right)\right)[u, u^{-1}]], b+uB\right)$
defined by $\tau_s\left(c\right) = \sum_{l \ge -\frac{\deg c -1}{2}}
u^{-l} \left(\tau_s\right)_l\left(c\right) $, where
$\left(\tau_s\right)_l\left(c\right) \in
C^k\left(C_0^{\ify}\left(P\right)\right)$, $k= \deg c +2l-1$, is
given by the formula:
\begin{multline}\label{deft}
\left(\tau_s\right)_l\left(c\right) \left(a_0, a_1, \ldots ,
a_{k}\right)\\= \sum \limits_{i=0}^{k}\left(-1\right)^i \left
\langle  c, \int \limits_{\De^{k+1}} \Tr_s a_0 e^{-t_0\Ab_s^2}
[\Ab_s, a_1] \dots
 \right.
\\
 \left.
e^{-t_i\Ab_s^2} \frac{d\Ab_s}{ds} e^{-t_{i+1}\Ab_s^2} \dots [\Ab_s,
a_k] e^{-t_{k}\Ab_s^2} dt_1\ldots dt_{k+1}  \right \rangle
\end{multline}
  Then we have the following:
\begin{lemma}\label{prtr}
\begin{equation}\label{trans}
\frac{d}{ds}\left(\theta_{\Ab_s}\right)\left(c\right)=\left(b+uB\right)\tau_s\left(c\right)
+ \tau_s\left(u\partial c\right)
\end{equation}
\end{lemma}
\begin{proof}
Consider the submersion $P\times [a,b] \to B\times [a,b]$ with
superconnection $d+\Ab_s$, $ s \in [a, b]$, where $d$ is de~Rham
differential on $[a,b]$. Here we view the interval $[a, b]$ as a
subset of $\Rb$ for the purpose of defining smooth functions, etc.
With the projection $p: B\times [a,b] \to B$ we can associate a
cochain $p^* \in \Hom^1 \left(\W_*\left(B\right), \W_*\left(B \times
[a, b] \right) \right)$, which is defined on currents of degree $k$
as $ \left(\int\right)^{t}$, where $\int: \W^*_c\left(B \times [a,
b]\right) \to \W^*\left(B \right)$ is the integration along the
fibers of $p$. Denote by $c_a$ a current on $B \times [a, b]$ whose
value on a form is a composition of restriction to $B\times a$ with
$c$, and similarly for $c_b$. Then if we denote by $u\partial$
differential in $\Hom^1 \left(\W_*\left(B\right), \W_*\left(B \times
[a, b]\right) \right)$ we have
\begin{equation}
\left(u \partial p^*\right) \left(c\right) =u \left(c_b -c_a\right).
\end{equation}
 Define also $e^* \in \Hom^0\left( C^*\left(C_0^{\ify}\left(P \times [a, b] \right)\right)[u,
u^{-1}]], C^*\left(C_0^{\ify}\left(P \right)\right)[u,
u^{-1}]]\right)$. $e^*$ is clearly a cocycle.
 The Proposition \ref{cocycle}
implies that $\theta_{\Ab_s+d}$ is a cocycle in $\Hom
\left(\W_*\left([0,1]\times B\right),
CC^*\left(C_0^{\ify}\left([0,1]\times P\right)\right) \right)$, and
hence
\begin{multline}
\left(b+uB+u\partial\right) \left(e^* \circ \theta_{d+\Ab_s} \circ
p^*\right)\left(c\right)= e^* \circ \theta_{d+\Ab_s} \circ
\left(u\partial p^*\right)\left(c\right)\\=u e^*\circ
\theta_{d+\Ab_s}\left(c_b-c_a\right)=u
\left(\theta_{\Ab_b}-\theta_{\Ab_a}\right)\left(c\right),
 \end{multline}
or
$\left(\theta_{\Ab_b}-\theta_{\Ab_a}\right)\left(c\right)=\left(b+uB+u\partial\right)
u^{-1}\left(e^* \circ \theta_{d+\Ab_s} \circ
p^*\right)\left(c\right)$. But it is easy to see that $
u^{-1}\left(e^* \circ \theta_{d+\Ab_s} \circ
p^*\right)\left(c\right)$ is exactly $\int \limits_a^b \tau_s
\left(c\right)ds$, and the statement of the Proposition follows.
\end{proof}

We will now study behavior of $\theta_{\Ab_s}$ and $\tau_s$ near
$s=0$.

\begin{lemma}\label{estimate}
Let $V_0$, $V_1$, \ldots $V_l$ are   operators acting on sections of
some vector bundle over  a manifold $M$. We assume that $V_i$ is a
composition of a pseudodifferential operator of order $v_i$ with a
compact support with a diffeomorphism of $M$, lifted to act on the
sections of the vector bundle, and $D$ is a first-order selfadjoint
pseudodifferential operator on $M$. If  $v_i=\max \{ \ord V_i, 0\}$
and $\sum v_i \le l$ then
\begin{equation}\label{vform}
\left| \int \limits_{\De^k} \Tr_s V_0 e^{-t_0s^2D^2} V_1 \ldots V_l
e^{-t_ls^2D^2} dt_1\ldots dt_l \right| =O\left(s^{-\sum v_i-\dim
M-1}\right)
\end{equation}
as $s \to 0$. If $V_i$ and $D$ depend continuously on some
parameters, the estimate is uniform on the compacts.
\end{lemma}
\begin{proof}
First, using the fact that  each $V_i$ is compactly supported
integral operator and exponential decay of the heat kernel off
diagonal we can replace $M$ by a compact manifold changing our
expression by at most $O\left(s^{\ify}\right)$.  If the case when
some of the operators $V_i$ have order $1$ and the other order $0$
this follows from \cite{gs89} and Weyl asymptotics. In general
replace each $V^i$ by $V_i'\left(1+D^2\right)^{v_i}$ with $V_i'$ of
order $0$. Then distribute powers of $\left(1+D^2\right)$ replacing
as needed $\left(1+D^2\right)^{\al}V_i'$ by
$\left(\left(1+D^2\right)^{\al}V_i'\left(1+D^2\right)^{-\al}\right)\left(1+D^2\right)^{\al}$
to reduce to the
  case when we have at most $\sum v_i+1$ operators $V_i$ of order 1 and the rest are of order $0$ .
\end{proof}
\begin{rem} More precise asymptotics $O\left(s^{-\sum \ord V_i-\dim M}\right)$
with no restrictions imposed in the previous lemma can be obtained
by methods of pseudodifferential calculus, compare  \cite{widom78,
widom80, gilkey95}.
\end{rem}
\begin{rem}\label{differ}
One can differentiate the expression in the left hand side of
\eqref{vform} by the parameters,obtaining again an expression of the
same kind, by the Duhamel's formula. By similar methods one can show
that the derivatives of order $\al$ are $O\left(s^{-\sum v_i-\dim
M-1-2\al}\right)$ as $s \to 0$.
\end{rem}
 \begin{prop} \label{asy} There exists a number $N$, depending on the
orders of components of superconnection as pseudodifferential
operators and $\dim B$ such that for all $l>N$ coefficients
$\theta_l$ for $u^{-l}$ in $\theta_{\Ab_s}$ have limit $0$ as $s \to
0$, and coefficients for $u^{-l}$ in $\tau_s$ are integrable near $s
= 0$.
\end{prop}
\begin{proof}
Let $\Ab^2=\sum \cf_{[i]}$, where $\cf_{[i]}$ is the component of
degree $i$. Then $\Ab_s^2=s^2D^2+ \sum_{i\ge1}
s^{\left(2-i\right)}\cf_{[i]}$.  Decomposing $[\Ab_s, a]$, $a \in
C_0^{\ify}\left(P\right)$ according to the form degree looks as
follows: $[\Ab_s, a]
=s\al_{[0]}\left(a\right)+\al_{[1]}\left(a\right)+ \ldots
s^{\left(1-i\right)}\al_{[i]}\left(a\right)+\ldots$ where
$\al_{[i]}\left(a\right)$ are pseudodifferential operators,
depending on $a$, and order of $\al_{[i]}\left(a\right)$ is $0$ for
$i=1$, $2$, and $k_i-1$ for $i\ge 2$. Using Duhamel's expansion we
can write component of cochain $\theta$ in $\Hom
\left(\W_p\left(B\right),
CC^k\left(C_0^{\ify}\left(P\right)\right)\right)$ as a sum of
finitely many terms of the form
\begin{equation*}
\int \limits_{\De^l} \Tr_s V_0 e^{-t_0D^2} V_1 \ldots V_l
e^{-t_lD^2} dt_1\ldots dt_l
\end{equation*}
 where each $V_i$ is either $\al_{[j]}$ for some
$j$ or $\cf_{[j]}$ for some $j$. The number of the terms of the form
$\al_{[j]}$ is $k$, and we denote these terms as $\al_{[i_1]},
\dots, \al_{[i_k]}$. Similarly denote the terms of the form
$\cf_{[j]}$ as $\cf_{[j_1]}, \dots \cf_{[j_m]}$, for some $m$. If we
rescale the superconnection, this term changes to
\begin{equation*}
\int \limits_{\De^l} s^{Q}\Tr_s V_0 e^{-t_0s^2D^2} V_1 \ldots V_l
e^{-t_ls^2D^2} dt_1\ldots dt_l
\end{equation*}
where
$Q=\left(1-i_1\right)+\ldots+\left(1-i_k\right)+\left(2-j_1\right)+\ldots+\left(2-j_m\right)=k-p+2m\ge
k-p$. Now notice that among the operators $V_i$ there are at most
$\dim B$ operators of nonzero form degree, and every operator of
zero form degree is bounded. Hence if $v_i$ is the order of $V_i$ as
a pseudodifferential operator, $\sum v_i \le v\dim B $, where
$v=\max { v_i}$. Since $l\ge k$  we see that if $k \ge v \dim B$ we
can apply the estimate of the Lemma \ref{estimate} and obtain that
our term is $O\left(s^{k-p-\left(\dim P -\dim B\right)-v\dim B
-1}\right)$, and if $k-p>n +v\dim B/+1$, the corresponding component
has a limit $0$ when $s \to 0$. Similarly one can show that with the
same bound on the degree $k-p$ we obtain $\tau_s$ is integrable at
$s=0$. For the future use notice that we can take $N=\dim
P+\left(v-1\right)\dim B+1$.
\end{proof}

\begin{prop} \label{deFF} Choose any even $k \ge N$. Define the cochain
$\Theta_{\Ab} \in \Hom^0 \left(\W_*\left(B\right), CC^* \left(
C_0^{\ify}\left(P\right) \right)\right)$ as $\sum u^{-l}
\left(\Theta_{\Ab}\right)_l$
\begin{equation}\label{dtet}
\left(\Theta_{\Ab}\right)_l=\begin{cases}\left(\theta_{\Ab}\right)_l & \text{ for } l< k\\
\left(\theta_{\Ab}\right)_k-\int \limits_0^1
\left(B+\partial\right)\left(\tau_s\right)_{k+1}ds
&\text{ for } l=k\\
0 &\text{ for } l>k
\end{cases}
\end{equation}
Then this cochain is a cocycle. Its cohomology class is independent
of $k$.
\end{prop}

\begin{proof}

Integrating equation \eqref{trans} and using the results of
Proposition \ref{asy} we obtain
\begin{equation}
\left(\theta_{\Ab}\right)_l=b \left(\delta \right)_{l-1}+
\left(B+\partial \right) \left(\delta_s \right)_{l+1}
\end{equation}
where $\delta_l= \int_0^1\left(\tau_s\right)_l ds$ and $l>k$. Set
also $\delta_l=0$ for $l < k$. Then
$\Theta_{\Ab}=\theta_{\Ab}-\left(uB+b+u\partial\right)\delta$ is a
cocycle, since $\theta_{\Ab}$ is one. If we change $k$ to $k'$,$k <
k'$ the cocycle $\Theta_{\Ab}$ will change by
$\left(b+uB+u\partial\right) \left(\sum_{k <i<k'} \delta_i\right)$
so its cohomology class remains the same.
\end{proof}

Assume now that we have a smooth family of superconnections
$\Ab(t)$, adapted to the family $D$. We then have the following
\begin{prop}\label{ho} Let $\Ab(t)$ be a continuous family of
superconnections adapted to the family $D$. Then
\begin{equation}
\Theta_{\Ab(1)} -\Theta_{\Ab(0)}= \left(b+uB+u\partial\right)T
\end{equation}

where $T_{\Ab(t)}  \in \Hom ^{-1}\left(\W_*\left(B\right),
CC^*\left(C_0^{\ify}\left(P\right)\right)\right)$ is defined by the
formula $T_{\Ab(t)}=u^{-1}\int \limits_{t \in [0,1]}e^* \circ
\Theta_{d+\Ab(t)} \circ p^*$. Here $p^* \in \Hom^1
\left(\W_*\left(B\right), \W_*\left(B\times [0,1]\right)\right)$ is
defined by
\begin{equation}
p^*\left(c\right) =  \left(\int
\limits_{B\times[0,1]/B}\right)^{t}c,
\end{equation}
$\Theta_{d+\Ab(t)}$ is a cocycle in $\Hom^0\left(\W_*\left(B \times
[0,1]\right), CC^*\left(C_0^{\ify}\left(P
\times[0,1]\right)\right)\right)$ and $e$ is the pull-back
$C_0^{\ify}\left(P\right) \to C_0^{\ify}\left([0,1]\times P\right)$.
Here   we use the same truncation to construct $\Theta_{d+\Ab(t)}$,
$\Theta_{\Ab(1)}$ and $\Theta_{\Ab(0)}$.
\end{prop}

\begin{proof}
 Notice that
the superconnection $d+\Ab(t)  $ satisfies all the conditions
necessary to construct the cocycle $\Theta_{d+\Ab(t) }$ by choosing
the appropriate truncation.  The rest of the  proof is the same as
the proof of the Lemma \ref{prtr}.
\end{proof}

To extend this construction to the case of nontrivial groupoid
action notice that this construction can be sheafified. Consider the
following sheaves on $B$: one is given by $U \mapsto
\W_*\left(U\right)$, the other is the sheafification of the presheaf
given by $U \mapsto CC^*\left(
C_0^{\ify}\left(\pi^{-1}\left(U\right)\right) \right)$. Both of
these sheaves are fine. The map $\Theta_{\Ab}$ constructed above
defines a morphism of these complexes of sheaves.

\subsection{}\label{groupoid}  We now review some notation and definitions
regarding groupoid algebras and cross-products by groupoids
\cite{haefliger79,cr99}. Let $F$ be  a soft sheaf on
$B=G^{\left(0\right)}$ with an action of $G$. This means that every
$\g \in G$ defines a map $F_{s\left(\g\right)} \to
F_{r\left(\g\right)}$, which is continuous. We denote by
$C^*\left(G, F\right)$ the complex of nonhomogeneous $G$-cochains
with values in $F$. It can be described as follows. Set for $n \ge
1$:
\begin{equation}
G^{\left(n\right)}=\{\left(\g_1, \g_2, \ldots, \g_n\right) \in G^n \
| \ s\left(\g_i\right)=r\left(\g_{i+1}\right), i=1, 2, \dots n-1\}
\end{equation}
We have for every $n$ a map $\ve_n : G^{\left(n\right)} \to B$
defined by
\begin{equation}
\ve_n \left(\g_1, \g_2, \ldots, \g_n\right) =r\left(\g_1\right)
\end{equation}
Set $C^n\left(G, F\right) =\G\left(G^{\left(n\right)};
\ve_n^*F\right)$. To define the coboundary operator we use the
simplicial maps $\de_i : G^{\left(n\right)} \to
G^{\left(n-1\right)}$. They are given for $n >1$ by the formula
\begin{equation}\label{sb}
\de_i \left(\g_1, \g_2, \ldots, \g_n\right)=\begin{cases}\left(\g_2,
\ldots, \g_n\right) &\text{ if }i=0
\\
\left(\g_1,\ldots,  \g_i\g_{i+1}\ldots, \g_n\right) &\text{ if } 1
\le i <
n-1 \\
\left(\g_1, \g_2, \ldots, \g_{n-1}\right) &\text{ if }i=n
 \end{cases} .\end{equation}
For $n=1$ $\de_0 \left(\g_1\right) =r\left(\g_1\right)$ and $\de_1
\left(\g_1\right)=s \left(\g_1\right)$. The coboundary $\de:
C^{n-1}\left(G,F\right) \to C^n\left(G, F\right)$ is given by
$\sum_{i=0}^n \left(-1\right)^i\de_i$ where the action of $G$ on $F$
is used to identify $\de_0^* \ve_{n-1}^* F$ with $\ve_n^*F$. We can
apply this construction if $F=\left(F^*,d\right)$ is a complex of
sheaves with a differential $d: F^* \to F^{*+1}$. In this case
$C^*\left(G, F^*\right)$ is a bicomplex;  we consider only finite
cochains in this bicomplex. To form a total complex of this
bicomplex we equip it with the differential $\pm \de + d$, where
$\pm=\left(-1\right)^n$  on $C^m\left(G, F^n\right)$.

If $A$ is a sheaf of algebras over $G$, one can form a $G$-sheaf
$CC^*\left(A\right)$, defined as sheafification of the presheaf $U
\mapsto CC^*\left(\G_c\left(U;A\right)\right)$. In our examples this
sheaf is always fine. Hence we can define a complex $C^*\left(G,
CC^*\left(A\right)\right)$. On the other hand one can form a
cross-product algebra $A\rtimes G$, defined as functions on $G$ with
$f\left(\g\right) \in A_{r\left(\g\right)}$. The product is given by
convolution taking the action of $G$ on $A$ into account. When $A$
is the sheaf of smooth functions on $B$ one obtains the convolution
algebra which we denote by $C_0^{\ify}(G)$. One can then define a
map of complexes
\begin{equation}\label{phi}
\Phi: C^*\left(G, CC^*\left(A\right)\right) \to CC^*\left(A \rtimes
G\right)
\end{equation}
This map has been constructed in \cite{cr99}, based on the
constructions in \cite{bn}, compare also \cite{ni90}, \cite{ft87},
\cite{gj93} . An important particular case is when $A$ is the sheaf
$C^{\ify}(B)$ of smooth functions on $B$. In this case there is a
canonical morphism of complexes of sheaves $\iota:
\W_*\left(B\right) \to CC^*\left( C_0^{\ify}\left(B\right)\right)$
defined on the current $c$ of degree $m$ by
\begin{equation}\label{iota}
\iota(c) (a_0, a_1, \ldots, a_m) =\frac{1}{m!}\left \langle c,
a_0da_1 \ldots da_m\right \rangle.
\end{equation}
Composing $\iota$ with the map in \eqref{phi} one obtains Connes'
map \cite{cn86, cn94}, also denoted by $\Phi$:
\begin{equation}
\Phi: C^*\left(G, \W_*\left(B\right)\right) \to
CC^*\left(C_0^{\ify}(G)\right).
\end{equation}

We will use the following properties of the map $\Phi$, compare
\cite{cr99}:
\begin{itemize}
\item If $A$, $B$ two $G$-algebras and $f:A \to B$ is a
$G$-homomorphism, we get a natural induced homomorphism,  which we
also denote by $f$, from $A \rtimes G$ to $B \rtimes G$. Also the
induced map $f^* :CC^*\left(B\right) \to CC^*\left(A\right)$ is
$G$-equivariant, and hence defines a map, also denoted by $f^*$,
from $C^*\left(G, CC^*\left(B\right)\right)$ to $C^*\left(G,
CC^*\left(A\right)\right)$. Then the following diagram commutes:
\begin{equation}
\xymatrix{
&C^*\left(G, CC^*\left(A\right)\right)  \ar[rr]^{\Phi}& &CC^*\left(A \rtimes G\right)\\
&&&\\&C^*\left(G, CC^*\left(B\right)\ar[uu]^{f^*}\right)
\ar[rr]^{\Phi}& &CC^*\left(B \rtimes G\right)\ar[uu]_{f^*}
         }
\end{equation}

 \item If $A$ is a filtered  algebra and the action of $G$
preserves filtration  then the complexes $C^*\left(G,
CC^*\left(A\right)\right)$ and $CC^*\left(A \rtimes G\right)$ also
are naturally filtered. The map $\Phi$ preserves this filtration.

\item Let $U_i$ be an open cover of $B$ and set $B'=\coprod U_i$.
Let $G'$ be the pull-back of the groupoid $G$ by the natural
projection $p:B' \to B$. Set also $A'=p^*A$. Then $A'$ is naturally
a $G'$-algebra. The cross-products $A\rtimes G$ and $A' \rtimes G'$
are naturally Morita equivalent. We also have a pull-back map
$C^*\left(G, CC^*\left(A\right)\right) \to C^*\left(G',
CC^*\left(A'\right)\right)$. Then the following diagram is
commutative up to homotopy:

\begin{equation}\label{morita}
\xymatrix{
&C^*\left(G', CC^*\left(A'\right)\right) \ar[rr]^{\Phi}& &CC^*\left(A' \rtimes G'\right)\\
&&&\\&C^*\left(G, CC^*\left(A\right)\right) \ar[rr]^{\Phi} \ar[uu]&
&CC^*\left(A \rtimes G\right)\ar[uu]
         }
\end{equation}
Here the vertical arrows are the isomorphisms induced by the
pull-back and Morita equivalence respectively.

\end{itemize}

\subsection{}\label{general}
 Let $G$ be an \'etale groupoid with the unit space
$G^{\left(0\right)} $. We say that a manifold $P$ is a $G$-space if
we are given submersion $\pi: P \to   G^{\left(0\right)}$,    and
for every $\g \in G$ we have a diffeomorphism $\g:
P_{r\left(\g\right)} \to P_{s\left(\g\right)}$. Here $P_b$, $b \in
B$ is a fiber $p^{-1}\left(b\right)$ over the point $b$. The
diffeomorphisms $\g$ should be compatible with the groupoid
structure. We assume also that the fibers are equipped with the
complete $G$-invariant Riemannian metric, and we are given a family
of $G$-invariant fiberwise Dirac operators acting on the sections of
$G$-equivariant bundle $\ce$. We define the map $\Theta_{\Ab}$
associated to the simplicial superconnection.

\begin{equation}
P^{\left(n\right)}=\{\left(p,\g_1, \g_2, \ldots, \g_n\right) \in
P\times G^n \ | \ \pi\left(p\right)=r\left(\g_1\right) \text{ and }
s\left(\g_i\right)=r\left(\g_{i+1}\right) \}
\end{equation}
For $n=0$ we set $P^{0}=P$.
 Note that $\pi$ induces a natural submersion $\pi_n: P^{\left(n\right)} \to
G^{\left(n\right)}$:
\begin{equation}
\pi_n \left(p,\g_1, \g_2, \ldots, \g_n\right) =\left(\g_1, \g_2,
\ldots, \g_n\right).
\end{equation}

For every $ 0\le i \le n$ we get a submersion map $\de_i:
P^{\left(n\right)} \to P^{\left(n-1\right)}$ defined by
\begin{equation}
\de_i \left(p,\g_1, \g_2, \ldots,
\g_n\right)=\begin{cases}\left(p\g_1, \g_2, \ldots, \g_n\right)
&\text{ if }i=0
\\
\left(p,\g_1,\ldots,  \g_i\g_{i+1}\ldots, \g_n\right) &\text{ if } 1
\le i <
n-1 \\
\left(p,\g_1, \g_2, \ldots, \g_{n-1}\right) &\text{ if }i=n
\end{cases} .\end{equation} The underlying map of the base
spaces is given by the formulas \eqref{sb}. We introduce also
submersion maps $\al_i: P^{\left(n\right)} \to P$, $i=0, 1, \ldots,
n$ defined by
\begin{equation}\label{defa}
\al_i\left(p,\g_1,\ldots, \g_n\right)=
\begin{cases}
p\g_1\ldots \g_i &\text{ if }i>0 \\
p                &\text{ if }i=0
\end{cases}.
\end{equation}
  The underlying maps of the base
spaces are $\al_i\left(\g_1,\ldots, \g_n\right)=s\left(\g_1\ldots
\g_i\right)$, with $\al_0\left(\g_1,\ldots,
\g_n\right)=r\left(\g_1\right)$. On each of the spaces
$P^{\left(n\right)}$ we consider the bundle
$\ce^{\left(n\right)}=\al_0^*\ce$. It is naturally isomorphic to
each of the bundles $\al_i^*\ce$, with the isomorphism given by the
action of $G$ on $\ce$.

We now give a definition of simplicial connection. Natural imbedding
of the simplex as a subset of $\Rb^n$ allows one to talk about
smooth functions, etc. on the simplex. By simplicial connection we
mean a collection of connections $\na^{\left(n\right)}$ on the $
\pi_n \times \id: P^{\left(n\right)} \times \De^{n} \to
B^{\left(n\right)}\times \De^n$ where $\De^n=\{\s_0, \ldots, \s_n
\in \Rb \ | \ \sum \s_i=1, \s_i \ge 0\}$ satisfying the following
compatibility conditions, compare \cite{dup78}:
\begin{equation}
 \left(\id \times \partial\right)^*\na^{\left(n\right)}   =\left(\de_i
\times \id\right)^*\na^{\left(n-1\right)}
\end{equation}
Here $\partial_i:\De^{n-1} \to \De^n$ is the  $i$-th face map given
by $\left(\s_0,\ldots \s_{n-1}\right)\mapsto \left(\s_0, \ldots
\s_{i-1}, 0, \s_i, \ldots \s_{n-1}\right)$. We also require that our
connection has the hollowing property:  in a local chart
$\na^{\left(n\right)}$ can be written as $d +\w$; we require that
$\w\left(\frac{\partial}{\partial \s_i}\right)=0$.

 Existence of
simplicial connections follows from the following construction. One
starts with an arbitrary connection $\na$ on the submersion $P \to
G^{\left(0\right)}$.
 Then
one can set
\begin{equation}
\na^{\left(n\right)}=\sum \limits_{i=0}^n \s_i \al_i^* \na  +d_{dR}
\end{equation}
where we set $\na^{\left(0\right)}=\na$. Here $d_{dR}$ is de Rham
differential on $\De^n$.

A $G$-invariant family $D$ naturally defines a family of fiberwise
operators on each of the submersion $\pi_n$, which we also denote by
$D$. We then define simplicial superconnection $\Ab$ on
$G$-submersion $\pi: P \to B$ as a collection of superconnections on
submersions $ \pi_n \times \id: P^{\left(n\right)} \times \De^{n}
\to G^{\left(n\right)}\times \De^n$, adapted to $D$ and satisfying
the compatibility conditions:
\begin{equation}\label{compat}
 \left(\id \times \partial_i\right)^*\Ab^{\left(n\right)}
=\left(\de_i \times \id\right)^*\Ab^{\left(n-1\right)}.
\end{equation}

We also require that they have the following properties. Locally the
superconnection $\Ab^{\left(n\right)}$ can be written as $d +\w$
where $\w$ is sum of differential forms on $G^{\left(n\right)}\times
\De^n$. We require that every component of $\w$ which has a positive
degree in $d\s$ variables also has a positive degree in
$G^{\left(n\right)}$ variables. It follows that
$\left(\Ab^{\left(n\right)}\right)^2$ also has this property. In
particular we see that there exists number $r>0$ such that every
component of $\w$ or of $\left(\Ab^{\left(n\right)}\right)^2$ of
degree $m$ in $d\s$ has degree at least $m/r$ in
$G^{\left(n\right)}$ direction. We also require that the degrees of
components of $\Ab^{\left(n\right)}$ as pseudodifferential operators
are bounded uniformly in $n$.

An example of such superconnection is given by the $\Ab_0=D+\na$,
with $\na$ an arbitrary simplicial connection. More precisely we set
$\Ab^{\left(n\right)}=D +\na^{\left(n\right)}$. Another example is
given by the Bismut simplicial superconnection described in the
section \ref{bismut}.

We now use simplicial superconnection $\Ab$ to construct the map
$\Phi_{\Ab}$. We start by constructing a cocycle $\{\Theta_{\Ab}^i\}
\in C^*\left(G, \Hom\left(\W_*\left(B\right),
CC^*\left(C_0^{\ify}\left(P\right)\right)\right)\right)$.

Introduce the cochains $p^n \in \Hom^n\left(\W_*\left(B\right),
\W_*\left(B\times \De^n\right)\right)$ by the formula
\begin{equation}
p^n\left(c\right) =  \left(\int \limits_{B\times \De^n/B}\right)^{t}
c.
\end{equation}
Here we view $\W_*\left(B \times \De^n\right)$ as a sheaf on $B$.
\begin{equation}
\theta^n_{\Ab}= u^{-n}\left(e_n\right)^* \circ
\theta_{\Ab^{\left(n\right)}}\circ p^n\in \Hom^{-n}
\left(\W_*\left(G^{\left(n\right)}\right),
C^*\left(C_0^{\ify}\left(P^{\left(n\right)}\right)\right)[u,
u^{-1}]], b+uB\right).
\end{equation}
Similarly one constructs
  $\tau_s^n$ by using the superconnection $\Ab^{\left(n\right)}$.
\begin{lemma} \label{vanish}There exists number $L$ such that
$\theta^n_{\Ab}=0$ and $\tau_s^n=0$ for $n >L$.
\end{lemma}
\begin{proof}

Let $U \in G^{\left(n\right)}$ be an open set, $c \in
\W_*\left(U\right)$, $a_i \in
C_0^{\ify}\left(\pi_n^{-1}\left(U\right)\right)$. Then

\begin{multline}
\left(e^* \circ
\theta_{\Ab^{\left(n\right)}}\right)_l\left(p^n\left(c\right)\right)
\left(a_0, a_1, \ldots , a_k\right)\\= \left \langle c,
\int\limits_{\De^n\times G^{\left(n\right)}|G^{\left(n\right)}} \int
\limits_{\De^k} \Tr_s a_0
e^{-t_0\left(\Ab^{\left(n\right)}\right)^2} [\Ab^{\left(n\right)},
a_1] \ldots [\Ab^{\left(n\right)}, a_k]
e^{-t_k\left(\Ab^{\left(n\right)}\right)^2} dt_1\ldots dt_k\right
\rangle.
\end{multline}

  Now due to the conditions imposed on
the superconnection every  component of  $[\Ab^{\left(n\right)},
a_i]$ of  degree $m$ in $d\s$ variables has degree at  least $m/k$
in the $G^{\left(n\right)}$ direction, and the same is true about
$\left(\Ab^{\left(n\right)}\right)^2$. Duhamel's expansion shows
that the same is true about
$e^{-t\left(\Ab^{\left(n\right)}\right)^2}$. Hence the component of
the expression under the integral which has degree $n$ in $d\s$
variables has degree at least $n/r$ in the $G^{\left(n\right)}$
direction. Since $\dim G^{\left(n\right)}= \dim B$ this expression
is $0$ for $n >r\dim B$. The same argument works for $\tau_s^n$.
\end{proof}
From this and the Proposition \ref{asy} we immediately obtain the
following:

 \begin{prop}   There exists a number $N$, depending on the
orders of components of superconnection as pseudodifferential
operators and $\dim B$ and independent of $n$ such that for all
$l>N$ $ \lim \limits_{s \to 0} \left(\theta_{\Ab_s}^n\right)_l =0 $
and $ \left(\tau_s^n\right)_l\left(c\right) $ is integrable near
$s=0$.
\end{prop}
We now can construct cochains $\Theta_{\Ab}^n \in \Hom^{-n}
\left(\W_*\left(G^{\left(n\right)}\right),
CC^*\left(C_0^{\ify}\left(P^{\left(n\right)}\right)\right)\right)$
as follows. Chose any $k
>N$, with $N$ as in the previous Proposition, and define the cochain
$\Theta_{\Ab^{\left(n\right)}}$ as in the Proposition \ref{deFF}. We
define then  $\Theta_{\Ab}^n=u^{-n}e^* \circ
\Theta_{\Ab^{\left(n\right)}} \circ p^n$. Notice that
$\Theta_{\Ab}^n =0$ for $n>r\dim B$. Hence we can view
$\Theta_{\Ab}^n$ as a cochain in $C^*\left(G,
\Hom\left(\W_*\left(B\right),
CC^*\left(C_0^{\ify}\left(P\right)\right)\right)\right)$. We then
have the following
\begin{lemma}\label{boundary}
\begin{equation}
\left(b+uB+u\partial\right)\Theta_{\Ab}^n=\left(-1\right)^{n-1}  \de
\Theta_{\Ab}^{n-1}.
\end{equation}
\end{lemma}
\begin{proof}
First, notice that as in the proof of Proposition \ref{deFF} we have
$e^* \circ \left(b+uB+u \partial\right)\Theta_{\Ab^{\left(n\right)}}
\circ p^n=0$. Hence $\left(b+uB+u\partial\right)\Theta_{\Ab}^n = e^*
\circ \Theta_{\Ab^{\left(n\right)}} \circ u\partial p^n$. We have
\begin{equation}
\partial p^n= \left(-1\right)^{n-1} \sum \limits_{i=0}^n \left(-1\right)^i\left(\id \times \partial_i\right)_* \circ
p^{n-1}.
\end{equation}
 The result of the Proposition then follows from this formula together with the compatibility conditions
 \eqref{compat}.
\end{proof}

\begin{thm}\label{mn}
The cochain $\{\Theta_{\Ab}^i \} = \sum_{n \ge 0} \Theta_{\Ab}^n \in
C^*\left(G, CC^*\left(C_0^{\ify}\left(P\right)\right)\right)$ is a
cocycle. The cohomology class of this cocycle is independent of the
choice of simplicial superconnection adapted to the family $D$.
\end{thm}
\begin{proof}
The first assertion -- that $\{\Theta_{\Ab}^i\}$ is a cocycle --
follows immediately from the Lemma \ref{boundary}. To see
independence of the connection let $\Ab$ be any  simplicial
superconnection and let $\Ab_0 =D+\na$, where
$\na=\{\na^{\left(n\right)}\}$ is an arbitrary simplicial
connection. Define then   $\Ab(t) =t\Ab +\left(1-t\right) \Ab_0$.
More precisely we set $\Ab(t)^{\left(n\right)}
=t\Ab^{\left(n\right)} +\left(1-t\right) \Ab_0^{\left(n\right)}$.
$\Ab(t)$ is then a simplicial superconnection for every $0\le t \le
1$, i.e. the equation \eqref{compat} is satisfied, as well as the
conditions listed after that equations. Moreover, we can use the
same number $r$ as defined there.
 It follows that we can use the construction of the Proposition
 \ref{ho} and construct for every $n$ a cochain $T_{\Ab(t)}^n=u^{-n}e^* \circ T_{\Ab(t)^{\left(n\right)}} \circ
 p^n$, which is $0$ for $n > r\dim B+1$. Hence we get a cochain $\{
 T_{\Ab(t)}^i\} \in C^*\left(G, \Hom \left(\W_*\left(B\right), CC^*\left(C_0^{\ify}\left(P\right)\right)\right)\right)$.
 As in the Lemma \ref{boundary} we obtain
 \begin{equation}\label{tr}
\left(b+uB+u\partial\right)T_{\Ab(t)}^n=\left(-1\right)^{n}  \de
T_{\Ab(t)}^{n-1}+\Theta_{\Ab}^n -\Theta_{\Ab_0}^n
\end{equation}
and hence $\left(b+uB+u\partial \pm
\de\right)\{T_{\Ab(t)}^i\}=\{\Theta_{\Ab}^i\}
-\{\Theta_{\Ab_0}^i\}$. We conclude that the cocycle
$\{\Theta_{\Ab}^i\}$ is cohomologous to the cocycle
$\{\Theta_{\Ab_0}^i\}$ for every simplicial superconnection $\Ab$.
\end{proof}

Now using the cup product
\begin{equation}
\cup: C^*\left(G, \Hom\left(\W_*,
CC^*\left(C_0^{\ify}\left(P\right)\right)\right)\right)\otimes
C^*\left(G, \W_*\right) \to C^*\left(G,
CC^*\left(C_0^{\ify}\left(P\right)\right)\right)
\end{equation}
we construct the map $C^*\left(G, \W_*\right) \to C^*\left(G,
CC^*\left(C_0^{\ify}\left(P\right)\right)\right)$ defined by $\alpha
\mapsto \{\Theta_{\Ab}^i\}\cup \alpha$. We obtain the map
$\Phi_{\Ab}$ by composing this map with the map $\Phi: C^*\left(G,
CC^*\left(C_0^{\ify}\left(P\right)\right)\right) \to
CC^*\left(C_0^{\ify}\left(P\right) \rtimes G\right)$. Explicitly we
define
\begin{equation}\label{defphi}
\Phi_{\Ab}\left(\al\right) =\Phi\left(\{\Theta_{\Ab}^i\}\cup
\alpha\right)
\end{equation}
The following theorem follows immediately from the Theorem \ref{mn}.
\begin{thm}
Let $\Ab$ be a simplicial superconnection adapted to the family $D$.
Then $\Phi_{\Ab}$ defined in the equation \eqref{defphi} is  a
cochain map of complexes. If $\Ab'$ is another simplicial
superconnection adapted to the family $D$, the maps $\Phi_{\Ab}$ and
$\Phi_{\Ab'}$ are   homotopic.
\end{thm}

\section{Fiberwise Pseudodifferential Operators and Bivariant Chern Character}\label{fpdo}
In this section we consider the properties of the algebra of the
proper fiberwise pseudodifferential operators which will be needed
in our discussion of the bivariant Chern character in the Section
\ref{biv}. We start by giving the general definitions in
\ref{define}. Then in \ref{deftr} we define the trace map on the
cyclic complex of this algebra. In \ref{connection} we discuss a
different construction of the trace map, involving connections.
Finally in \ref{equality} we show that the two maps are the same up
to homotopy.

\subsection{}\label{define} Let $\cf$ be a $G$-equivariant bundle over $P$. We
assume that $\cf$ is $\Zb_2$-graded with the grading given by the
operator $\g\in \End\left(\cf\right)$. In this section we consider
the algebra $\Psi\left(\cf\right)$ of the fiberwise
pseudodifferential operators on the submersion $P \to B$ of order
$0$ acting on the sections of the bundle $\cf$ which are even with
respect to the grading and whose Schwartz kernel is compactly
supported. The $\Zb_2$ grading we use here is  induced by the
grading of $\cf$.

This algebra also has a natural filtration by the order of
pseudodifferential operators. This filtration induces corresponding
filtration on the cyclic homology complex of the algebra
$\Psi\left(\cf\right)$. We denote by
$F^{-k}CC_*\left(\Psi\left(\cf\right)\right)$ the subcomplex of the
cyclic complex $CC_*\left(\Psi^0\right)$ generated by the $
\bigoplus \limits_{j_0+\ldots j_l \le -k} \Psi^{j_0}\left(\cf\right)
\otimes \ldots \otimes \Psi^{j_l}\left(\cf\right)$. Thus we obtain
filtration
\[
CC_*\left(\Psi\left(\cf\right)\right)=F^0CC_*
\left(\Psi\left(\cf\right)\right)\supset
F^{-1}CC_*\left(\Psi\left(\cf\right)\right) \supset
F^{-2}CC_*\left(\Psi\left(\cf\right)\right) \supset \ldots
\]

It follows from Goodwillie's theorem \cite{goo85}, compare also
\cite{cq95a}, that the inclusion
$F^{-i}CC_*\left(\Psi\left(\cf\right)\right) \to
F^{-1}CC_*\left(\Psi\left(\cf\right)\right)$ is a quasiisomorphism
for every $i$.

We will also consider the dual complexes
$F_kCC^*\left(\Psi\left(\cf\right)\right)$, so that we have
\[
CC^*\left(\Psi\left(\cf\right)\right)=F_{0}CC^*\left(\Psi\left(\cf\right)\right)\subset
F_{1}CC^*\left(\Psi\left(\cf\right)\right) \subset
F_{2}CC^*\left(\Psi\left(\cf\right)\right) \subset \ldots
\]
Groupoid $G$ acts on $\Psi\left(\cf\right)$, preserving the order
filtration. We can then form the cross-product algebra $
\Psi\left(\cf\right) \rtimes G$, which inherits filtration from the
algebra $\Psi\left(\cf\right)$. We use this filtration to construct
the complexes $F^{-k}CC_*\left(\Psi\left(\cf\right) \rtimes
G\right)$ and $F_{k}CC^*\left(\Psi\left(\cf\right) \rtimes
G\right)$. Note also that the map $\Phi$ defines a map of complexes:
\begin{equation}
\Phi: C^*\left(G, F_iCC^*\left(\Psi\left(\cf\right)\right)\right)
\to F_iCC^*\left(\Psi\left(\cf\right)\rtimes G\right).
\end{equation}

\subsection{} \label{deftr}

Now \emph{ assuming that $P \to B$ is a fibration} we construct for
every  $n >  \dim\left(P\right)-\dim\left(B\right)$ a map
\begin{equation}
\tau: F^{-n}CC_* \left(\Psi\left(\cf\right)\rtimes G\right) \to CC_*
\left(C_0^{\ify}\left(G\right)\right).
\end{equation}
 The construction is as follows, compare \cite{crmo01}. Consider a covering $U_i$ of $B$
trivializing the fibration. Set $B'=\coprod U_i$. Also set
$P'=\coprod \pi^{-1}U_i$. The pull-back groupoid  $G'$ then acts on
the trivial fibration $\pi': P' \to B'$.  Morita equivalence induces
an isomorphism $ CC_* \left(\Psi\left(\cf\right)\rtimes G\right) \to
CC_* \left(\Psi\left(\cf'\right)\rtimes G'\right) $, where $\cf'$ is
the pull-back of $\cf$ to $P'$. This isomorphism preserves
filtrations. This allows us to assume that the fibration is trivial,
$P=F\times B$ with $\pi$ given by the projection on the second
factor. Notice that since the fibration is trivial the action of $G$
on $B$ induces an action on the fibration by $\left(f, b\right)
\g=\left(f, b\g\right)$. We will call this action the product
action. This action is different in general from the action arising
from the original action on $P$.

We will need to consider a bigger algebra $\Psi\De= \Psi
\De\left(\cf\right)=\Psi \left(\cf\right) \rtimes \Diff$ instead of
the algebra $\Psi\left(\cf\right)$, where $\Diff$ is the group of
smooth families of the fiberwise diffeomorphisms, viewed as a
discrete group. An element of this algebra can be viewed as a
fiberwise integral operator on the fibration with a distributional
kernel. We will denote by $\Psi \De \rtimes_t G$ the cross product
of the algebra $\Psi \De$ by the product action. The algebra $\Psi
\De $ contains a dense subalgebra defined as an algebraic tensor
product of the algebra $K$ of the smoothing integral operators on
the fiber with the smooth compactly supported functions on $B$. The
corresponding cross-product algebra is the algebraic tensor product
$C_0^{\ify}\left(G\right) \otimes K$, and we define $\tau_0$ by

\begin{equation}
\tau_0 \left(\left(a_0\otimes k_0\right)\otimes \ldots \left(a_i
\otimes k_i\right)\right) =\Tr_s\left(k_0\ldots k_i\right)
a_0\otimes a_1\ldots a_i
\end{equation}
 where $a_i \in C_0^{\ify}\left(G\right)$, $k_i \in K$.   The map $\tau_0$ then has a
 unique extension to the cyclic complex $F^{-n}CC_*(\Psi \De \rtimes_t G)$.

 We now consider the general case, when the action on the
 fibration is not necessarily the product. In this case the action
 of $\g \in G$ is given by
 \begin{equation}\label{twac}
 \left(f, b\right)\g=\left(\s\left(\g\right) f, b\g\right)
\end{equation}
where $\s\left(\g\right)$ is a diffeomorphism of the fiber $F$.
Notice that $\s$ has to   satisfy the cocycle condition:
\begin{equation}
\s\left(\g_1 \g_2\right)=\s\left(\g_2\right) \s\left(\g_1\right)
\end{equation}

We will denote by  $\Psi \De \rtimes G$ the cross-product of $\Psi
\De$ by $G$ with the action given by \eqref{twac}.  These two
cross-products are isomorphic as filtered algebras. Explicitly the
isomorphism $I:\Psi\De \rtimes  G \to \Psi \De \rtimes_t G$  is
given by
\begin{equation}\label{act}
\left(I\left(f\right)\right)\left(\g\right)=f\left(\g\right)\s\left(\g\right)
\end{equation}

We can now define the trace map  $\tau$ by
\begin{equation} \label{deftau}
\tau= \tau_0 \circ I_*\circ i_*
\end{equation}
where $i:\Psi\left(\cf\right)\rtimes G \to \Psi\De \rtimes G$ is the
inclusion map.

\subsection{}\label{connection} On the other hand a choice of simplicial connection $\na$
as in the section \ref{construction} provides us with the map
\begin{equation}
\ct_{\nabla}: C^*\left(G, \W_*\left(B\right)\right) \to
F_{k}CC^*\left(\Psi\right)
\end{equation}
where $k$ depends on the choice of the superconnection $\na$, or
more precisely on order of it as a vertical pseudodifferential
operator. The definition is parallel to the definition of the map
$\Theta_{\Ab}$ but only the connection part is used instead of the
full superconnection. Since the definition is very close to the
construction in the section \ref{construction}, we give just a brief
outline here.

 Notice that  $\Psi\left(\cf\right)$ naturally can be
viewed as compactly supported sections of  a sheaf over $B$
$U\mapsto  \Psi^U\left(\cf\right)$, where $\Psi^U\left(\cf\right)$
is the algebra of proper fiberwise pseudodifferential operators on
the submersion $\pi^{-1}\left(U\right) \to U$. We can then form a
sheaf $CC^*\left(\Psi\left(\cf\right)\right)$ as the sheafification
of the presheaf $U\mapsto
CC^*\left(\left(\Psi^U\left(\cf\right)\right)_c\right)$.

 First we consider the case of the trivial
action of $G$. In that case for any choice of connection $\na$ on $P
\to B$ we can define a cocycle $\cl_{\na} \in \Hom
\left(\W_*\left(B\right),
F_{m}CC^*\left(\Psi\left(\cf\right)\right)\right)$ by the formula
similar to the equation \eqref{defj} and $m$ is specified below in
\eqref{mbo}. Explicitly we define
\begin{equation}
\cl_{\na}: c \mapsto \sum \limits_{l \ge  -\frac{\deg c}{2} }
u^{-l}\cl_l\left(c\right)
\end{equation}
where $\cl_l\left(c\right) \in
C^k\left(\Psi\left(\cf\right)\right)$, $k= \deg c +2l$, is given by
the formula:
\begin{multline}
\cl_l\left(c\right) \left(a_0, a_1, \ldots , a_k\right)\\= \left
\langle c, \int \limits_{\De^k} \Tr_s a_0 e^{-t_0\na^2} [\na, a_1]
\ldots [\na, a_k] e^{-t_k\na^2} dt_1\ldots dt_k\right \rangle
\end{multline}

Here we view $\nabla^2$ as a $2$-form on $B$ with values in the
fiberwise differential operators and define $e^{-t \nabla^2}$ by the
usual series.

If $\nabla$ locally looks like $d +\w$ where $\w$ is a one-form on
$B$ with values in the proper fiberwise pseudodifferential operators
of order $v\ge 1$. Then $\na^2$ is a vertical operator of order at
most $2v$. If $a_0\otimes\ldots \otimes a_k$ above is in
$F^{-m}CC_*\left(\Psi\left(\cf\right)\right)$ then the expression
under the integral above is a compactly supported fiberwise
pseudodifferential operator of order at most $v \dim B-m$.  Hence
the expression under the integral will be well-defined if
\begin{equation}\label{mbo}
m> \left(v-1\right)\dim B+\dim P.
\end{equation}

Next given a  simplicial connection $\na^{\left(n\right)}$ we
construct a cocycle $\{\cl_{\na}^i\}\in C^*\left(G, \Hom
\left(\W_*\left(B\right),
F_{m}CC^*\left(\Psi\left(\cf\right)\right)\right)\right)$.   The
$n$-th component $\cl_{\na}^n \in \ve_n^*\Hom^{-n}
\left(\W_*\left(G\right), F_{m}CC^*\left(\Psi\left(\cf\right)
\right)\right)$ is given by $\cl_{\na}^n=u^{-n}e^*\circ
\cl_{\na^{\left(n\right)}} \circ p^n$. We assume that order of
$\na^{\left(n\right)}$ as a pseudodifferential operator is bounded
by $v$, independent from $n$. Then for $\cl_{\na}^n$  to be
well-defined it is sufficient to have $m
>\left(v-1\right)\left(\dim B+n\right)+\dim P$. However, as before we have $\cl_{\na}^n=0$ for $n>\dim B$, by exactly
the same argument as in the section \ref{construction}. This implies
that we indeed get a cochain in $C^*\left(G, \Hom
\left(\W_*\left(B\right),
F_{m}CC^*\left(\Psi\left(\cf\right)\right)\right)\right)$, and also
that any $m>2\left(v-1\right)\dim B +\dim P$ can be used for all
$n$. The same argument shows that this cochain is indeed a cocycle.
We can now construct a map of complexes $C^*\left(G,
\W_*\left(B\right)\right) \to CC^*\left(G,
F_{m}CC^*\left(\Psi\left(\cf\right)\right)\right)$ by taking the cup
product with $\{\cl_{\na}^i\}$. Following it by the map $\Phi$ we
obtain our map $\ct_{\na}$. It is a map of complexes and different
choice of simplicial connection leads to a chain-homotopic map.
Explicitly we have
\begin{equation}
\ct_{\na}\left(\al\right)=\Phi\left(\{\cl_{\na}^i\}\cup \al\right).
\end{equation}

\subsection{}\label{equality}
  In this section we
prove the following result:
\begin{thm} \label{trthm}The following diagram is commutative up to homotopy.
\begin{equation}\label{trace}
\xymatrix{
&& &F_mCC^*\left(\Psi\left(\cf\right) \rtimes G\right)\\
&&&\\&  C^*\left(G, \W_*\left(B\right)\right)
\ar[rruu]^{\ct_{\nabla} \ } \ar[rr]^{\Phi \ \ \ \ \ \ } & & CC^*
\left( C_0^{\ify}\left(G\right) \right)\ar[uu]_{\left(
\tau\right)^t}
         }
\end{equation}
\end{thm}
\begin{proof}

First we reduce this statement to  the case of trivial fibration $P
\to B$. Let $P'$, $B'$ and $G'$ be as in section \ref{deftr}, and
$\Psi'$ -corresponding sheaf of algebras of pseudodifferential
operators on $P'\to B'$. We then have for every $n$ a natural
projection $\rho: \left(G'\right)^{\left(n\right)} \to
G^{\left(n\right)}$, and $\left(P'\right)^{\left(n\right)}$ is a
pull-back, as a bundle, of $P^{\left(n\right)}$ under this
projection. We can construct a simplicial connection $\na'$ by
setting
$\left(\na'\right)^{\left(n\right)}=\rho^*\na^{\left(n\right)}$. It
is easy to see then that $\{\cl_{\na'}^i\}=\rho^* \{\cl_{\na}^i\}$.
This implies that the diagram
\begin{equation}
\xymatrix{
&C^*\left(G', \W_*\left(B'\right)\right) \ar[rr]^{\cdot \cup \{\cl_{\na'}^i\} \ \ \ }& &C^*\left(G', F_mCC^*\left(\Psi\left(\cf'\right)\right)\right)\\
&&&\\&C^*\left(G, \W_*\left(B\right)\right)
\ar[uu]^{\rho^*}\ar[rr]^{\cdot \cup \{\cl_{\na}^i\}\ \ \ }&
&C^*\left(G', F_mCC^*\left(\Psi\left(\cf\right)\right)\right)
\ar[uu]_{\rho^*}
         }
\end{equation}
commutes. This, together with \eqref{morita} implies commutativity
of the diagram
\begin{equation}
\xymatrix{
&C^*\left(G', \W_*\left(B'\right)\right) \ar[rr]^{\ct_{\na'}}& &F_mCC^*\left(\Psi\left(\cf'\right) \rtimes G'\right)\\
&&&\\&C^*\left(G, \W_*\left(B\right)\right)\ar[uu]^{\rho^*}
\ar[rr]^{\ct_{\na}}& &F_mCC^*\left(\Psi\left(\cf\right) \rtimes
G\right)\ar[uu]
         }
\end{equation}
 up to homotopy. Here the right vertical arrow is induced by the Morita equivalence.
This, together with the definition of the map $\tau$ and
\eqref{morita} implies that the statement in the general case
follows from the statement in the case of the trivial fibration.

Next consider the case of the trivial fibration with the product
action.  We can chose a trivial connection  on  fibrations
$P^{\left(n\right)}\to G^{\left(n\right)}$, i.e. the one given by
the de Rham differential $d$ with respect to the decomposition
$P^{\left(n\right)}=F\times G^{\left(n\right)}$. This clearly
defines a simplicial connection which we denote $\na_0$.
Commutativity of the diagram \eqref{trace} in this case is clear.

Consider now the general case: $P=F \times B$, the action is given
by the formula \eqref{act}, and let $\na_1$ be an arbitrary
simplicial connection. The statement will follow from the previous
observations and the following result: the diagram

\begin{equation}
\xymatrix{
& & &F_mCC^*\left(\Psi\De \rtimes G\right)\\
&C^*\left(G, \W_*\left(B\right)\right)\ar[urr]^{\ct_{\na_1}}
\ar[drr]^{\ct_{\na_0}}&&\\& & &F_mCC^*\left(\Psi \De \rtimes_t
G\right) \ar[uu]_{I^*}
         }
\end{equation}
is commutative up to homotopy. To prove this consider the algebra
$M_2\left(\Psi \De\right)$. The product action of $G$ on $\Psi \De$
induces an action on the algebra $M_2\left(\Psi \De\right)$ which we
denote as $m \mapsto m^{\g}$. Consider now another action of $G$ on
this algebra:
\begin{equation}
\g\left(m\right) =  \begin{bmatrix}1 &0 \\ 0 &\s\left(\g\right)
\end{bmatrix} m^{\g}\begin{bmatrix}1 &0 \\ 0 &\s\left(\g\right)
\end{bmatrix}^{-1}
\end{equation}
 We use this action to form the cross-product algebra $M_2\left(\Psi
\De\right)\rtimes G$. Notice that we have homomorphisms $i_0: \Psi
\De \rtimes_t G \to M_2\left(\Psi \De\right)\rtimes G$ and $i_1:
\Psi \De \rtimes G \to M_2\left(\Psi \De\right)\rtimes G$ defined by
\begin{align}
i_0\left(f\right)\left(\g\right)&= \begin{bmatrix} f\left(\g\right) &0 \\
0 &0 \end{bmatrix}\\
i_1\left(f\right)\left(\g\right)&= \begin{bmatrix} 0 &0 \\
0 &f\left(\g\right) \end{bmatrix}
\end{align}
Notice that these homomorphisms are induced by the $G$-homomorphisms
of the corresponding algebras. Let $\Lambda$ be an automorphism of
$M_2\left(\Psi \De\right)~\rtimes~G$ given by
\begin{equation}
\Lambda\left(f\right) = \begin{bmatrix}0 &-1\\1 &0 \end{bmatrix}f
\begin{bmatrix}0 &-1\\1 &0 \end{bmatrix}^{-1}.
\end{equation}
 Notice that
\begin{equation}\label{inst}
i_0\circ I =\Lambda \circ i_1
\end{equation}
  We can construct a
new simplicial connection $\nabla_0 \oplus \nabla_1$ on
$\pi_*\left(\cf\right) \oplus \pi_*\left(\cf\right)$. Its
commutation with elements of $M_2\left(\Psi \De\right)$ is given by
\begin{equation}
\left[ \nabla_0\oplus \nabla_1, \begin{bmatrix} a&b \\ c &d
\end{bmatrix} \right]=\begin{bmatrix} [\na_0, a] &[\na_0,b]-b\left(\na_1-\na_0\right)\\
[\na_0,c]+\left(\na_1-\na_0\right)c & [\na_1, d]\end{bmatrix}
\end{equation}
and the curvature is
 \begin{equation}
 \left(\na_0\oplus \na_1\right)^2= \begin{bmatrix} \na_0^2 &0\\
 0 & \na_1^2 \end{bmatrix}
 \end{equation}
One way to view this construction is to replace $P$ with $P\sqcup
P$, fibered over $B$, with the  $G$ action on the first copy being
the product action and on the second copy being the one given by
equation \eqref{act}.  The algebra $M_2\left(\Psi \De\right)$ is
then the cross-product algebra for this $G$-fibration.

 We can use  connection $\na_0\oplus \na_1$ to construct the map
\begin{equation}
\ct_{\na_0\oplus \na_1}: C_*\left(G, \W_*\left(B\right)\right) \to
F_mCC^*\left(M_2\left(\Psi \De\right) \rtimes G\right)
\end{equation}
for $m$ large enough. We then have
\begin{align}
\left(i_0\right)^*\ct_{\na_0\oplus \na_1}&=\ct_{\na_0}\\
\left(i_1\right)^*\ct_{\na_0\oplus \na_1}&=\ct_{\na_1}
\end{align}
Here we have used the fact that the map $\Phi$ is compatible with
homomorphisms. It follows these equalities and equation \eqref{inst}
that
\begin{equation}
I^*\circ \ct_{\na_0}=I^*\circ\left(i_0\right)^* \circ
\ct_{\na_0\oplus \na_1}=\left(i_1\right)^*\circ \Lambda^*\circ
\ct_{\na_0\oplus \na_1}
\end{equation}
We now notice that $\Lambda^* \circ \ct_{\na_0\oplus \na_1}$ is
homotopic to $ \ct_{\na_0\oplus \na_1}$. Indeed, let $\Lambda_t$ be
an automorphism of $M_2\left(\Psi \De\right)\rtimes G$ given by
 \begin{equation} \Lambda_t\left(f\right)  = \begin{bmatrix}\cos t &-\sin t\\\sin t &\cos t \end{bmatrix}f
\begin{bmatrix}\cos t &\sin t\\ -\sin t &\cos t \end{bmatrix}
\end{equation}
Then $\Lambda_0 =\id$, $\Lambda_{\frac{\pi}{2}}=\Lambda$. Since we
have $\frac{d}{dt}\Lambda_t\left(f\right) =[\lb,
\Lb_t\left(a\right)]$, where $\lb=
\begin{bmatrix}
                        0& -1&\\
                        1&  0&
                        \end{bmatrix}
                 $
these maps satisfy $\frac{d}{dt}\Lb_t^*   = \Lb_t^*L_{\lb} $, where
$L_{\lb}$ is  the action of the derivation $\ad_{\lb}$ on the cyclic
complex. Notice that this derivation preserves the order filtration.
We then have a homotopy formula
\begin{equation}
L_{\lb}=[\left(B+b\right), H_{\lb}]
\end{equation}
for the action of $L_{\lb}$, \cite{cn85, goo85, ge93a}. Explicit
formulas for $H_{\lb}$ are not important here, we just need to know
that it preserves the order filtration, since $L_{\lb}$ does. Since
$\ct_{\na_0\oplus \na_1}$ is a map of complexes we conclude that
 \begin{multline}
 \frac{d}{dt}\Lb_t^*\circ \ct_{\na_0\oplus
\na_1}\left(c\right)=\\\left(b+uB\right)\left(H_{\lb}\Lb_t^*\ct_{\na_0\oplus
\na_1}\right)\left(c\right)-\left(H_{\lb}\Lb_t^*\ct_{\na_0\oplus
\na_1}\right)\left(\pm \de+u\partial\right).
 \end{multline}
 Integrating from
$0$ to $\pi/2$ provides us with desired homotopy and finishes the
proof.
\end{proof}

\section{Comparison with the bivariant character}\label{biv}
We start by reviewing in \ref{bch} V.~Nistor's construction of the
bivariant Chern character. Then in \ref{quih} we explain the
construction of the quasihomomorphism $\psi_D$. Finally in
\ref{lower} we show the commutativity of the lower triangle in
\eqref{comm'} up to homotopy.

\subsection{}\label{bch} In \cite{ni91, ni93a} V. Nistor constructed a bivariant Chern
character of a quasihomomorphism. We will need the following from
this construction. Given algebras $A$ and $B$, where $B$ is equipped
with filtration $B=B_0 \supset B_1 \ldots$, with $B_iB_j \subset
B_{i+j}$.   We can use filtration as before to introduce, cyclic
complexes $F^iCC_*\left(B\right)$ and $F_iCC^*\left(B\right)$.

 A quasihomomorphism $\Upsilon$ is a pair of
homomorphisms $\upsilon_0$, $\upsilon_1: A \to B$ such that
 \begin{equation}
 \upsilon_1\left(a\right)-\upsilon_0\left(a\right) \in B_{1}
 \end{equation}
 V. Nistor constructed a sequence of
maps $c_i\left(\Upsilon\right): CC_*\left(A\right) \to
CC_{*}^{-i}(B)$, $i=1, 2, \ldots$, which satisfy the following
properties:
\begin{itemize}
\item The first map $c_1\left(\Upsilon\right)$ is given by the formula
\begin{equation}\label{defc1}
  c_1\left(\Upsilon\right)= \frac{1}{2}\left(\left(\upsilon_1\right)_* -\left(\upsilon_0\right)_*\right)
\end{equation}
We introduce the coefficient $\frac{1}{2}$ since we consider below
the symmetrized version of quasihomomorphism associated to an
operator, compare \cite{cn85}.
 \item
If $r_i:F^{-i}CC_*\left(B\right)\to F^{-1} CC_*\left(B\right)$ is
the natural inclusion, then $r_i\circ c_i\left(\Upsilon\right)$   is
homotopic to $c_1\left(\Upsilon\right)$ for $i \ge 1$.
\end{itemize}

The construction proceeds as follows. Starting with the algebra $A$
one constructs certain canonical algebra $QA$ with filtration,
together with canonical quasihomomorphism $j=\left(j_0, j_1\right)$
from $A$ to $QA$. One constructs canonical maps $s_i:
F^{-1}CC_*\left(QA\right) \to F^{-i}CC_*\left(QA\right)$ as follows:
$s_1=\id$, $s_{i+1}=s_i+[b+uB, H_i]$, where $H_i$ is a certain
canonical endomorphism of $ CC_*\left(QA\right)$, which preserves
filtration. We note that $s_i$ defines a homotopy equivalence of the
complexes $F^{-1}CC_*\left(QA\right) $ and
$F^{-i}CC_*\left(QA\right)$ for every $i$, with the homotopy inverse
  given by the inclusion $F^{-i}CC_*\left(QA\right) \to F^{-1}CC_*\left(QA\right)$.
 One then can define the bivariant Chern character of the canonical quasihomomorphism
$j$. This is a sequence of  canonical maps $c_k\left(j\right):
CC_*\left(A\right) \to F^kCC_*\left(QA\right)$ where
$c_1\left(j\right)=\frac{1}{2}\left(\left(i_1\right)_*-\left(i_0\right)_*\right)$
and $c_k\left(j\right)=s_k\circ c_1$ for $k \ge 1$. Now the
quasihomomorphism $\Upsilon$ from $A$ to $B$ defines canonically a
homomorphism $h:QA \to B$, which preserves filtration. The bivariant
Chern character of $\Upsilon$ then is defined as
$c_i\left(\Upsilon\right)=h_* \circ c_i\left(j\right)$.

 We introduce also $c^i\left(\Upsilon\right)=c_i\left(\Upsilon\right)^t$.
 Assume now
that $A$ and $B$ are $G$-algebras, filtration on $B$ is
$G$-invariant and $\upsilon_0$, $\upsilon_1$ are $G$-equivariant.
 In
this case these maps have the following naturality property. A
$G$-equivariant quasihomomorphism induces a quasihomomorphism from
$A\rtimes G$ to $B \rtimes G$, which we also denote by $\Upsilon$.
Since all the steps in the construction of the bivariant Chern
character are canonical   $c_i\left(\Upsilon\right)$ can be chosen
to be $G$-equivariant. Hence its transposed induces a map
$C^*\left(G, F_iCC^*\left(B\right)\right)  \to C^*\left(G,
CC^*\left(A\right)\right)$.
\begin{prop}\label{natl}
The following diagram is commutative up to homotopy for every $i\ge
1$:
\begin{equation}\label{natd}
\xymatrix{
&C^*\left(G, CC^*\left(A\right)\right)\ar[rr]^{\Phi}   & &CC^*\left(A \rtimes G\right) \\
&&&\\&C^*\left(G,
F_iCC^*\left(B\right)\right)\ar[uu]^{c^i\left(\Upsilon\right)}
\ar[rr]^{\Phi} & &F_iCC^*\left(B\rtimes
G\right)\ar[uu]^{c^i\left(\Upsilon\right)}
  }
\end{equation}
\end{prop}
\begin{proof}
We start by showing the commutativity of the diagram
\begin{equation}
\xymatrix{
&C^*\left(G, CC^*\left(A\right)\right)\ar[rr]^{\Phi}    & &CC^*\left(A \rtimes G\right) \\
&&&\\&C^*\left(G,
F_iCC^*\left(QA\right)\right)\ar[uu]^{c^i\left(j\right)}\ar[rr]^{\Phi}
& &F_iCC^*\left(QA\rtimes G\right)\ar[uu]^{c^i\left(j\right)}
  }
\end{equation}
First notice that since $QA$ is constructed canonically it is indeed
a $G$-algebra, so the statement makes sense. For $i=1$ the
commutativity is clear.  Hence it is enough to show the
commutativity up to homotopy of the diagram
\begin{equation} \label{qdiag}
\xymatrix{
&C^*\left(G, F_1CC^*\left(QA\right)\right)\ar[rr]^{\Phi}    & &F_1CC^*\left(QA \rtimes G\right) \\
&&&\\&C^*\left(G,
F_iCC^*\left(QA\right)\right)\ar[uu]^{s^i}\ar[rr]^{\Phi} &
&F_iCC^*\left(QA\rtimes G\right)\ar[uu]^{s^i}
  }
\end{equation}

Here both vertical arrows are transposed of the maps $s_i$ mentioned
above. The right vertical arrow is constructed for the algebra
$QA\rtimes G$. The left is constructed  from the transposed of $s_i$
for the  algebra $QA$. Since the construction is canonical it is
$G$-equivariant, and hence defines a map $C^*\left(G,
F_iCC^*\left(QA\right)\right) \to C^*\left(G,
F_1CC^*\left(QA\right)\right)$. Notice that both of these maps are
homotopy equivalences of complexes, with the inverses induced by
inclusion. This statement follows from the Goodwillie's theorem
\cite{goo85} for the right vertical arrow. For the left arrow we
notice that the homotopy $H_i$ above is constructed canonically and
hence is $G$-equivariant. As a result $H_i$ defines a
filtration-preserving endomorphism of $C^*\left(G,
CC^*\left(QA\right)\right)$. This implies that $s^i$ is indeed a
homotopy equivalence. Note that since we consider only the finite
cochains in the complex $C^*\left(G, CC^*\left(QA\right)\right)$ the
fact that $s^i: F_iCC^*\left(QA\right) \to F_1CC^*\left(QA\right)$
is a homotopy equivalence does not imply that
$C^*\left(G,F_iCC^*\left(QA\right)\right) \to C^*\left(G,
F_1CC^*\left(QA\right)\right)$ is a homotopy equivalence; we need
the more precise version of the argument given above.

Hence the commutativity of the diagram \eqref{qdiag} follows from
the commutativity of the diagram
\begin{equation}
\xymatrix{
&C^*\left(G, F_1CC^*\left(QA\right)\right)\ar[dd]\ar[rr]^{\Phi}    & &F_1CC^*\left(QA \rtimes G\right)\ar[dd] \\
&&&\\&C^*\left(G, F_iCC^*\left(QA\right)\right) \ar[rr]^{\Phi}  &
&F_iCC^*\left(QA\rtimes G\right)
  }
\end{equation}
where the vertical arrows are induced by the transposes of the
inclusions. This is true since $\Phi$ preserves filtrations.

Now to deduce the general case we note that the homomorphism $QA \to
B$ being canonical is $G$-equivariant. Hence the commutativity of
the diagram \eqref{natd} follows from the compatibility of $\Phi$
with homomorphisms.

\end{proof}

\subsection{}\label{quih}
For our bundle $\ce$ consider the algebra $\Psi\left(\ce \oplus
\ce\right)$ of even operators on $\ce \oplus \ce$. Here
$\Zb_2$-grading is defined as follows. The bundle $\ce$ has a
grading given by the operator $\g\in \End\left(\ce\right)$. This
induces the grading on $\ce\oplus \ce$, and hence on
$\Psi\left(\ce\oplus \ce\right)$,   given by the operator
\begin{equation}
\G=\begin{bmatrix} \g &0 \\
0&-\g\end{bmatrix}
 \end{equation}
We will write elements of $\Psi\left(\ce \oplus \ce\right)$ as
$2\times 2$ matrices of pseudodifferential operators on $\ce$.
Consider now the algebra $\Psi_r$  defined as a subalgebra of
$\Psi\left(\ce \oplus \ce\right)$ whose principal zero-order symbol
has a form $\begin{bmatrix} a&0\\0 &0
\end{bmatrix}$, where $a$ is a scalar endomorphism of $\ce$.
 We associate with
the family $D$ a $G$-equivariant quasihomomorphism from the algebra
$C_0^{\ify}\left(P\right)$ to the  algebra $\Psi_r$. Here we use the
filtration on $\Psi_r$ induced by the filtration on $\Psi\left(\ce
\oplus \ce\right)$.

We now define quasihomomorphism $\psi_D$   by the following
formulas:
\begin{align}
\psi_0\left(a\right) &=U_D\begin{bmatrix} a &0 \\ 0&0\end{bmatrix}U_D^{-1}\\
\psi_1\left(a\right)&=\begin{bmatrix} a &0 \\ 0&0\end{bmatrix}
\end{align}
Here $U_D$ is constructed as follows. Let $Q$ be a family of proper
pseudodifferential operators forming a parametrix of $D$. One can
always chose $Q$ to be $G$-equivariant, as follows for example from
explicit formulas for parametrix \cite{mw94a}. Then $S_0=1-QD$ and
$S_1=1-DQ$ are proper smoothing $G$-equivariant fiberwise operators.
Define then $U_D$ by the formula
 \begin{equation}
 U_D=\begin{bmatrix} D &S_1 \\ S_0&-\left(1+S_0\right)Q\end{bmatrix}
 \end{equation}
Its inverse $U_D^{-1}$ is given by the following explicit formula:
 \begin{equation}
  U_D^{-1}=\begin{bmatrix} \left(1+S_0\right)Q&S_0 \\ S_1&-D\end{bmatrix}
 \end{equation}
It is easy to see that $\psi_1\left(a\right)-\psi_0\left(a\right)
\in \Psi_r^{-1}$. Moreover $\psi_0$ and $\psi_1$ are
$G$-equivariant, since $D$ is $G$-equivariant.
 Notice also that $U_D$ is odd with respect to
$\G$:
 \begin{equation}
 U_D\G=-\G U_D
 \end{equation}
and in particular $\psi_0(a)$ is an even operator.

Since the quasihomomorphism $\psi_D$ is $G$-equivariant, it also
defines a quasihomomorphism, also denoted by $\psi_D$, from
$C_0^{\ify}(P) \rtimes G$ to $\Psi_r \rtimes G$.  One now gets for
every $i$ the map of complexes $c_i\left(\psi_D\right):
CC_*\left(C_0^{\ify}\left(P\right)\right) \to
F^{-i}CC_*\left(\Psi_r\rtimes G\right)$, satisfying the above listed
properties.  We now define the bivariant Chern character
$\Ch\left(D\right): CC_*\left(C_0^{\ify}\left(P\right)\rtimes G
\right) \to CC_*\left(C_0^{\ify}\left(G\right)\right)$ by the
formula
\begin{equation}
\Ch\left(D\right) =\tau \circ c_i\left(\psi_D\right)
\end{equation}
where $\tau$ is defined in \eqref{deftau} and $i >\dim P -\dim B$.
It is clear that different choices of $i$ give homotopic maps.
\begin{rem} Write $D=\begin{bmatrix} 0 &D^{+}\\ D^{-}
&0\end{bmatrix}$ , where the decomposition is with respect to the
grading $\g$. We can use above formulas applied to $D^{+}$, omitting
however the $\frac{1}{2}$ factor from \eqref{defc1} and using
ordinary trace instead of supertrace,  to construct the bivariant
Chern character $\Ch\left(D^+\right)$. Similarly we can construct
$\Ch\left(D^-\right)$. It is easy to see that
\begin{equation}
\Ch\left(D\right)=\frac{1}{2}\left(\Ch\left(D^+\right)-\Ch\left(D^-\right)\right)
\end{equation}
It is however easy to see that $\Ch\left(D^+\right)$ is homotopic to
$-\Ch\left(D^-\right)$, and hence $\Ch\left(D\right)$ and
$\Ch\left(D^+\right)$ are homotopic.
\end{rem}
 Different choices of $Q$ are homotopic, and hence lead to
the homotopic maps of complexes. For our purposes it will be
convenient to chose $Q$ so that it commutes with $D$. That this is
possible again follows  from the explicit construction in
\cite{mw94a}. In this case all the entries in $U_D$ commute with $D$
and hence $U_D$ commutes with $\begin{bmatrix} D &0\\ 0 &D
\end{bmatrix}$.

\subsection{}\label{lower}

We start by proving the following:
\begin{prop}
Let $\Ab$ be any simplicial superconnection adapted to the family
$D$, and let $\na$ be any simplicial connection on $\ce$, even with
respect to $\g$. Then the following diagram is commutative up to
homotopy.
\begin{equation}
\xymatrix{
&& &CC^*\left(C_0^{\ify}\left(P\right)\rtimes G\right)\\
&&&\\&  C^*\left(G, \W_*\left(B\right)\right) \ar[rruu]^{
\Phi_{\Ab}}\ar[rr]^{\ct_{\widetilde{\na}} } & &
F^{m}CC^*\left(\Psi_r \rtimes G\right)\ar[uu]_{
c^m\left(\psi_D\right)}
         }
\end{equation}

\end{prop}
\begin{proof}
First notice that by the Theorem \ref{mn} we can assume  that
$\Ab=D+\na$.  Consider now the superconnection $\widetilde{\Ab}=\Ab
\oplus \Ab$ on the bundle $\widetilde{\ce}=\ce \oplus \ce$. By this
we mean that we lift naturally $\widetilde{\ce}$ to
$P^{\left(n\right)}$ and define a simplicial superconnection
$\widetilde{\Ab}$ by
$\widetilde{\Ab}^{\left(n\right)}=\Ab^{\left(n\right)}\oplus\Ab^{\left(n\right)}$.
It is adapted to the $G$-invariant family $\widetilde{D} =D \oplus
D$. We can now construct  a cocycle $\{X_{\Ab}^i\} \in C^*\left(G,
\Hom \left(\W_*\left(B\right),
F_1CC^*\left(\Psi_r\right)\right)\right)$ by exactly the same
formulas as $\{\Theta_{\Ab}^i\}$. Since for $A \in \Psi_r$ the
commutator $[\widetilde{D}, A]$ has order $0$, all the estimates
used for the construction of $\{\Theta_{\Ab}^i\}$ still hold, and
the formulas make sense. We use  for $\{X_{\widetilde{\Ab}}^i\}$ the
same truncation as for $\{\Theta_{\Ab}^i\}$. Proof of the Theorem
\ref{mn} applies here as well and shows that if one replaces
$\widetilde{\Ab}$ by another simplicial superconnection adapted to
$\widetilde{D}$ one obtains a cohomologous cocycle.

 Recall that $c^1\left(\psi_D\right)$ defines a
$G$-equivariant map $F_1CC^*\left(\Psi_r\right) \to
CC^*\left(C_0^{\ify}\left(P\right)\right)$, and hence a map from
$C^*\left(G, F_1CC^*\left(\Psi_r\right)\right)$ to
$C^*\left(G,CC^*\left(C_0^{\ify}\left(P\right)\right)\right)$. We
then have the following:
\begin{lemma}
The cocycles $c^1\left(\psi_D\right) \circ
\{X_{\widetilde{\Ab}}^i\}$ and $\{\Theta_{\Ab}^i\}$ are
cohomologous.
\end{lemma}
\begin{proof}
First notice that both $\psi_0^* \circ \{X_{\widetilde{\Ab}}^i\}$
and $\psi_1^* \circ \{X_{\widetilde{\Ab}}^i\}$ are well defined, and
not just their difference. It is clear that $\psi_1^* \circ
\{X_{\widetilde{\Ab}}^i\} =\Phi_{\Ab}$. On the other hand  the
identity
\begin{multline}
\Tr_s U_D a_0 U_D^{-1} e^{-t_0\Ab^2} [\Ab, U_Da_1U_D^{-1}] \ldots
[\Ab, U_Da_kU_D^{-1}] e^{-t_k\Ab^2}=\\ -\Tr_s a_0
e^{-t_0\left(U_D^{-1} \Ab U_D\right)^2} [\left(U_D^{-1} \Ab
U_D\right), a_1] \ldots [\left(U_D^{-1} \Ab U_D\right), a_k]
e^{-t_k\left(U_D^{-1} \Ab U_D\right)^2}
\end{multline}
together with a similar identity for components of $\tau_s$ imply
that $\psi_0^* \circ \{X_{\widetilde{\Ab}}^i\}=-\psi_1^* \circ
\{X_{U_D^{-1}\widetilde{\Ab} U_D}^i\} $. Here
$U_D^{-1}\widetilde{\Ab} U_D$ is the simplicial superconnection
defined by $\left(U_D^{-1}\widetilde{\Ab}
U_D\right)^{\left(n\right)}=
U_D^{-1}\left(\widetilde{\Ab}\right)^{\left(n\right)} U_D$. It is
also adapted to the operator $\widetilde{D}$. Explicitly we have
 \begin{equation}
 U_D^{-1}\widetilde{\Ab}U_D=\begin{bmatrix} D &0\\0 &D
 \end{bmatrix}+U_D^{-1}\begin{bmatrix} \na &0\\0 &\na
 \end{bmatrix} U_D
 \end{equation}
Since the cohomology class of $\{X_{\widetilde{\Ab}}^i\}$ is
independent of the superconnection we obtain that
$\psi_0^*\{X_{\widetilde{\Ab}}^i\}$ is cohomologous to
$-\psi_1^*\{X_{\widetilde{\Ab}}^i\}$, and the statement of the Lemma
follows.
\end{proof}
Consider now a map $X:\al \mapsto
\Phi\left(\{X_{\widetilde{\Ab}}^i\} \cup \alpha\right)$. It follows
that the composition of $X$ with the map $c^1\left(\psi_D\right):
F_1CC^*\left(\Psi_r~\rtimes~G\right) \to
CC^*\left(C_0^{\ify}\left(P\right)\rtimes G\right)$ is homotopic to
$\Phi_{\Ab}$.
  Let $r^m:
F_1CC^*\left(\Psi_r \rtimes G\right) \to F_mCC^*\left(\Psi_r \rtimes
G\right)$ be the be the transposed of the inclusion $r_m$.   Since
$c^m\left(\psi_D\right) \circ r^m$ is homotopic to
$c^1\left(\psi_D\right)$, to complete the proof it is sufficient to
show that $r^m \circ X$ is homotopic to $\ct_{\widetilde{\na}}$. To
establish this we consider a homotopy
$\widetilde{\Ab}_t=t\widetilde{\Ab}
+\left(1-t\right)\widetilde{\na}=t\widetilde{D} +\widetilde{\na}$
between $\widetilde{\Ab}$ and $\widetilde{\na}$. The existence of
this homotopy clearly implies that $r^m\circ
\{X_{\widetilde{\Ab}}^i\}$ and $\{\cl_{\widetilde{\na}}^i\}$ are
cohomologous after one establishes an analogue of the identity
\eqref{tr} in this context. Explicitly, one needs to show that the
analogue of the cochain $T$ is well define. This is done exactly as
before, with the only difference being that the estimate of the
Lemma \ref{estimate} has to be replaced by the following:
\begin{lemma}
Let $V_0$, $V_1$, \ldots $V_l$ be pseudodifferential operators with
a compact support acting on sections of some vector bundle over a
manifold $M$. Let $D$ be a first-order selfadjoint
pseudodifferential operator on $M$.  Assume that $\sum \ord V_i
<-\dim M$. Then the expression
\begin{equation}
\left| \int \limits_{\De^k} \Tr_s V_0 e^{-t_0s^2D^2} V_1 \ldots V_l
e^{-t_ls^2D^2} dt_1\ldots dt_l \right|
\end{equation}
is bounded uniformly in $s$.
\end{lemma}
\begin{proof}
Let $\phi \in C_0^{\ify}\left(M\right)$ be such that $\phi V_0
=V_0$. Denote $ V_0 e^{-t_0s^2D^2} V_1 \ldots V_l e^{-t_ls^2D^2}$ by
$A$. Then $\Tr_s A=\Tr_s A \phi$. If $v \sum \ord V_i$, then the
operator $A\phi D^{-v}$ is bounded uniformly in $s$. Let now $\psi
\in C_0^{\ify}\left(M\right)$ be such that $\phi \psi =\phi$. Then
$|\Tr_s A|= |\Tr_sA \phi D^{-v}D^v\phi| \le \|A \phi\| |\Tr
D^v\phi|$, which is bounded uniformly.
\end{proof}
This estimate implies that for sufficiently small $m$ $T$ is well
defined. This completes the proof of the Proposition.
\end{proof}

We now can prove the main result of this section:
\begin{thm}\label{ltr}
The diagram
\begin{equation}
\xymatrix{
&& & CC^*\left(C_0^{\ify}\left(P\right)\rtimes G\right)\\
&&&\\&  C^*\left(G, \W_*\left(B\right)\right) \ar[rruu]^{\quad
\Phi_{\Ab}}\ar[rr]^{ \Phi} & &
CC^*\left(C_0^{\ify}\left(G\right)\right)
\ar[uu]_{\Ch\left(D\right)^t}
         }
\end{equation}
is commutative up to homotopy.
\end{thm}
\begin{proof}
The Proposition \ref{natl} together with the previous Proposition
immediately imply commutativity of the diagram
\begin{equation}
\xymatrix{
&& & CC^*\left(C_0^{\ify}\left(P\right)\rtimes G\right)\\
&&&\\&  C^*\left(G, \W_*\left(B\right)\right) \ar[rruu]^{\quad
\Phi_{\Ab}}\ar[rr]^{ \ct_{\widetilde{\na}}} & & F_mCC^*\left(\Psi_r
\rtimes G\right) \ar[uu]_{c^m\left(\psi_D\right)}
         }
\end{equation}
Combined with the Theorem \ref{trthm} this gives the proof of our
Theorem.
\end{proof}

\section{Bismut Superconnection and Short-time limit}\label{bismut}
In this section we obtain a topological expression for $\Phi_{\Ab}$,
thus finishing the proof of the Theorem \ref{mainthm}. In order to
do this we construct simplicial Bismut superconnection $\Ab$ for
which $\Phi_{\Ab_s}$, with $\Ab_s$ the rescaled superconnection, has
a limit when $s \to 0$.

First we will need the following result:
\begin{prop}\label{homin}
Let $\Ab$ be a simplicial superconnection. Consider the family of
maps $\Phi_{\Ab_s}$, where $\Ab_s$ is the rescaled superconnection.
Then the maps obtained for different values of $s>0$ are homotopic.
\end{prop}
\begin{proof}
Consider first the case of the trivial groupoid action. In this case
we have in the notations of the Proposition \ref{deFF}
\begin{equation}
\frac{d}{ds}
\Theta_{\Ab_s}\left(c\right)=\left(b+uB+u\partial\right)\sum
\limits_{ l \ge -\frac{\deg c -1}{2}}^k u^{-l}
\left(\tau_s\right)_l\left(c\right)
\end{equation}
where $\tau_s$ is defined in the equation \eqref{deft}. Denote the
cochain appearing in the right hand side by $T_{\Ab_s}$. In the
general case construct the cochain $\{T_{\Ab_s}^i\}\in C^*\left(G,
\Hom \W_*\left(B\right),
CC^*\left(C_0^{\ify}\left(P\right)\right)\right)$ by
$T_{\Ab_s}^n=u^{-n} e^* \circ T_{\Ab^{\left(n\right)}_s} \circ p^n$.
Calculations as in the Lemma \ref{boundary} then show that
\begin{equation}
\frac{d}{ds}\{\Theta_{\Ab_s}^i\}=\left(b+uB+u\partial \pm \de\right)
\{T_{\Ab_s}^i\}.
\end{equation}
Integration of this identity shows that cocycles
$\{\Theta_{\Ab_s}^i\}$ obtained for different values of $s$ are
cohomologous. This implies that corresponding maps $\Phi_{\Ab_s}$
are homotopic.
\end{proof}

We now proceed to construct the simplicial Bismut superconnection.
Recall that we are given a $G$-submersion $\pi: P \to B$.  The
fibers are equipped with the complete $G$-invariant Riemannian
metric, and we are given a family of $G$-invariant fiberwise Dirac
operators acting on the sections of $G$-equivariant Clifford module
bundle $\ce$. We start by reviewing the case of the trivial groupoid
action, see \cite{bgv} for the details. In this case one needs to
make the following choices. First one needs to choose the horizontal
distribution   on the submersion $P \to B$. By this we mean choice
of a smooth subbundle $H$ of the tangent bundle $TP$ such that
$TP\cong H \oplus TP/B$. Here $TP/B$ is the vertical tangent bundle.
The choice of $H$ is, of course not unique. Note that at each point
$p \in P$ the set of all possible $H_p$ has a natural structure of
an affine space based on the vector space $ \Hom \left(
T_{\pi\left(p\right)}B, T_pP/B\right)$. Hence the set of all
horizontal distributions thus has a natural structure of an affine
space. A choice of a horizontal distribution $H$ is equivalent to
the choice of projection $P_H: TP \to TP/B$. We have the following
relation:
\begin{equation}
P_{tH_1+\left(1-t\right)H_2}=tP_{H_1}+\left(1-t\right)P_{H_2}
\end{equation}

  Choice of $H$ allows one also
to construct a canonical connection $
\na^{P/B}=\na^{P/B}\left(H\right)$ on $TP/B$, whose restriction on
the fibers of $\pi$ coincides with Levi-Civita connection. This
connection can be defined by the relation
\begin{multline}\label{npb}
\left(\na^{P/B}_XY,Z\right)=\frac{1}{2}\left(\left(P_H[X,Y],
Z\right)-\left([Y, Z], P_HX\right)+ \right.\\ \left. \left(P_H[Z,
X], Y\right) +X\left(Y, Z\right)-Z\left(P_HX, Y\right)+Y\left(Z,
P_HX\right)\right)
\end{multline}

From this equation it is clear that the correspondence $H\mapsto
\na^{P/B}\left(H\right)$ is an affine map.  One then needs to chose
a connection $\na^{\ce}$ on the bundle $\ce$ which is compatible
with the connection $\na^{P/B}\left(H\right)$. Notice that here
again the vertical part of this connection is determined uniquely.
Now using the horizontal distribution $H$ and a connection
$\na^{\ce}$ one can construct a connection $\nabla=\nabla^H$ on
$\pi_* \ce$ \cite{bgv}.   The following facts are easy consequences
of the definitions:
\begin{itemize}
\item If $H_1$, $H_2$ are two horizontal distributions and
$\na^{\ce}_1$ and $\na^{\ce}_2$ are two connections on $\ce$
compatible with the connections $\na^{P/B}\left(H_1\right)$ and
$\na^{P/B}\left(H_2\right)$ respectively. Then the connection
$t\na^{\ce}_1+\left(1-t\right)\na^{\ce}_2$ is compatible with
$\na^{P/B}\left(tH_1+\left(1-t\right)H_2\right)$.

\item Let $\na^{H_i}$ be constructed using connection $\na^{\ce}_i$,
$i=1, 2$. Then the connection $\na^{tH_1+\left(1-t\right)H_2}$
constructed using $t\na^{\ce}_1 + \left(1-t\right)\na^{\ce}_2$ is
equal to $t\na^{H_1}+\left(1-t\right)\na^{H_2}$.
\end{itemize}

The Bismut superconnection is then defined as
\begin{equation}
\Ab =D +\na^H-\frac{1}{4}c\left(T^H\right)
\end{equation}
where $c$ denotes the Clifford multiplication by a vertical vector
field and $T$ is a $2$-form on $B$ with the values in the vertical
vector fields  given by the curvature of the distribution $H$.
Explicitly it is described as follows. Let $X$, $Y$ be two vector
fields on $B$, and $X^H$, $Y^H$ their lifts to the horizontal vector
fields on $P$. Then
\begin{equation}
T^H\left(X, Y\right)=-\left([X^H, Y^H]-[X, Y]^H\right)
\end{equation}

We now explain how one can extend this construction to the
simplicial context. We fix the horizontal distribution $H$. Fix also
a connection $\na^{\ce}$ compatible with $\na^{P/B}$. We will now
use this data to construct a simplicial superconnection such that
each component $\Ab^{\left(n\right)}$ is a Bismut superconnection on
the corresponding submersion. We start by constructing a horizontal
distribution on each of the submersions $P^{\left(n\right)} \times
\De^n \to G^{\left(n\right)}\times \De^n$.  Recall the maps $\al_i:
P^{\left(n\right)} \to P$ defined in \eqref{defa}. Define
$\widetilde{\al}_i$ as the composition of the projection
$P^{\left(n\right)}\times \De^n \to P^{\left(n\right)}$ with the map
$\al_i$ Let $\s_i$ be the barycentric coordinates on $\De^n$; below
we view them as functions on $P^{\left(n\right)} \times \De^n$.
Define the distribution $H^{\left(n\right)}$ as
\begin{equation}
H^{\left(n\right)}=\sum \limits_{i=0}^n \s_i \left(d
\widetilde{\al}_i\right)^{-1}\left(H\right)
\end{equation}
The corresponding projection  is given by the formula
\begin{equation}
P_{H^{\left(n\right)}}=\sum \limits_{i=0}^n \s_i
\left(\left(\widetilde{\al}_i\right)^{*}\left(P_H\right)\oplus
0\right).
\end{equation}
Notice that this horizontal distribution satisfies the compatibility
conditions
\begin{equation}
 \left(\id \times \partial_i\right)^{-1}H^{\left(n\right)}
=\left(\de_i \times \id\right)^{-1}H^{\left(n-1\right)}.
\end{equation}
Using the formula \eqref{npb} and the horizontal distribution
$H^{\left(n\right)}$ we can construct for every $n$ the canonical
connection $\left(\na^{P/B}\right)^{\left(n\right)}$ on the vertical
bundle for the submersion $P^{\left(n\right)}\times \De^n \to
G^{\left(n\right)}\times \De^n$. It follows from the formula
\eqref{npb} that these connections are given explicitly by the
formula
\begin{equation}
\left(\na^{P/B}\right)^{\left(n\right)}=\sum \limits_{i=0}^n \s_i
\al_i^*\na^{P/B}\left(H\right) +d=\sum \limits_{i=0}^n \s_i
\left(\al_i^*\na^{P/B}\left(H\right) +d\right)
\end{equation}
where $d$ is  de Rham differential on $\De^n$. It is clear that
these connections satisfy the simplicial compatibility conditions.
It follows then that the connections on
$\left(\na^{\ce}\right)^{\left(n\right)}$ on $
\widetilde{\al}_0^*\ce$ defined by the formula
\begin{equation}
\left(\na^{\ce}\right)^{\left(n\right)} = \sum \limits_{i=0}^n \s_i
\al_i^*\na^{\ce} +d
\end{equation}
are compatible with the connections
$\left(\na^{P/B}\right)^{\left(n\right)}$. We can now construct for
every $n$ Bismut superconnection $\Ab^{\left(n\right)}$ on the
submersion $P^{(n)} \times \De^n \to G^{(n)} \times \De^n$ using the
distribution $H^{\left(n\right)}$ and connection
$\left(\na^{\ce}\right)^{\left(n\right)} $. The simplicial
compatibility conditions \eqref{compat} are clearly satisfied.

 Also, the only terms involving $d\s$ in the superconnection form are
the ones arising from $T^{H^{\left(n\right)}}$. But it is easy to
see that $T^{H^{\left(n\right)}}\left(\frac{\partial}{\partial
\s_i}, \frac{\partial}{\partial \s_i}\right)=0$, and hence there are
no terms of positive degree in $d \s$ and degree $0$ in
$G^{\left(n\right)}$ direction. It follows that the conditions
stated after the equation \eqref{compat} are satisfied as well; $r$
can be taken to be $1$.

To state the next proposition we   follow conventions and notations
from \cite{bgv}. Let $R^{\left(n\right)}$ be the curvature of the
connection $\left(\na^{P/B}\right)^{\left(n\right)}$. One can then
define the $\widehat{A}\left(R^{\left(n\right)}\right)$, a
differential form on $P^{(n)}$, where
$\widehat{A}\left(x\right)=\det^{\frac{1}{2}}\left(\frac{x/2}{\sinh(x/2)}
\right)$. The structure of the Clifford module on $\ce$ allows one
to view $R^{\left(n\right)}$ as a $2$-form on $P$ with values in the
endomorphisms of $\ce$. The twisting curvature
$F^{\left(n\right)}_{\ce/\cs}$ of
$\left(\na^{\ce}\right)^{\left(n\right)}$ is defined as
$F^{\left(n\right)}_{\ce/\cs}=\left(\left(\na^{\ce}\right)^{\left(n\right)}\right)^2-R^{\left(n\right)}$

We now have the following
\begin{prop} \label{limit} a) The cocycle $\Theta_{\Ab_s}$ has a limit when $s \to
0$ given by
\begin{equation}
\lim \limits_{s \to 0} \Theta_{\Ab_s}^n = \sum \limits_{l}u^{-n-l}
\phi_l^n
\end{equation}
where $\phi_l^n \in \Hom^{-n}
\left(\W_*\left(G^{\left(n\right)}\right),
CC^*\left(C_0^{\ify}\left(P^{\left(n\right)}\right)\right)\right)$
is defined by the formula
\begin{multline}
\phi_l^n\left(c\right) \left(a_0, a_1, \ldots, a_m\right) =\\
\frac{\left(2 \pi i\right)^{-\frac{\dim P-\dim B}{2}}}{m!}\left
\langle  \int \limits_{\left(P^{\left(n\right)}\times
\De^n\right)/P^{\left(n\right)}}
\widehat{A}\left(R^{\left(n\right)}\right) \Ch
\left(F^{\left(n\right)}_{\ce/\cs}\right)\pi_n^*c, a_0da_1 \ldots
da_m \right \rangle
\end{multline}
$m=\deg c+ 2l +n$.

  \noindent b) The cochain $T_{\Ab_s}$ defined in the
proof of the Proposition \ref{homin} has a limit when $s \to 0$.
\end{prop}
\begin{proof}
From the construction of $\Theta_{\Ab_s}^n$ it is clear that it is
enough to compute $\lim_{s\to 0} \Theta_{\Ab_s}$ with
$\Theta_{\Ab_s}$ defined in the Proposition \ref{deFF}. From the
equation \eqref{dtet} we obtain
\begin{equation}
\left(\Theta_{\Ab_s}\right)_l=\begin{cases}\left(\theta_{\Ab_s}\right)_l
& \text{ for } l< k\\ \left(\theta_{\Ab_s}\right)_k-\int \limits_0^s
\left(B+\partial\right)\left(\tau_t\right)_{k+1}dt
&\text{ for } l=k\\
0 &\text{ for } l>k
\end{cases}
\end{equation}
We see that it is enough to study the behavior of
$\left(\theta_{\Ab_s}\right)_l$ when $s \to 0$. Now the standard
application of the Getzler's calculus \cite{bgv} shows that
\begin{multline}
\lim \limits_{s\to 0} \left(\theta_{\Ab_s}\right)_l\left(c\right)\left(a_0, a_1,\ldots, a_m\right)=\\
  \frac{\left(2 \pi i\right)^{-\frac{\dim P-\dim B}{2}}}{m!}\left
\langle \widehat{A}\left(R\right) \Ch
\left(F_{\ce/\cs}\right)\pi^*c, a_0da_1 \ldots da_m \right \rangle
\end{multline}
where $R$ is the curvature of the connection $\na^{P/B}$, and
$F_{\ce/\cs}$ is the twisting curvature of $\na^{\ce}$. The
statement of the part a) follows. The proof of the part b) is
analogous to the proof of the local regularity of $\eta$-forms
\cite{bf86}.  Namely one shows that $T_{\Ab_s}$ has asymptotic
expansion in powers of $s$, starting with $s^{-1}$. Then one shows
that the coefficient for the leading term is $0$.
\end{proof}

Note that the above Proposition implies that we can represent
$\lim_{s \to 0} \Theta_{\Ab_s}$ as a composition of two maps.

The first one is $\widetilde{\pi}^* : \W_*(G^{(n)}) \to
\W_*(P^{(n)})$ defined by
\begin{equation}
\widetilde{\pi}^*:c \mapsto \left(2 \pi i\right)^{-\frac{\dim P-\dim
B}{2}}\int \limits_{\left(P^{\left(n\right)}\times
\De^n\right)/P^{\left(n\right)}}
\widehat{A}\left(u^{-1}R^{\left(n\right)}\right) \Ch
\left(u^{-1}F^{\left(n\right)}_{\ce/\cs}\right) \wedge \pi_n^*(c)
\end{equation}
It is easy to see that this map commutes with the differential $\pm
\de$ and hence defines a map of complexes $C^*\left(G,
\W_*\left(B\right)\right) \to C^*\left(G, \W_*\left(P\right)\right)$
To identify this map on the level of cohomology notice that the
collection of forms $ \widehat{A}\left(R^{\left(n\right)}\right) \Ch
\left(F^{\left(n\right)}_{\ce/\cs}\right)$  on $P^{(n)} \times
\De^n$ defines a closed simplicial form in the sense of
\cite{dup78}. The cohomology class of this form is the product of
$\widehat{A}_G\left(TP/B\right)$, the equivariant
$\widehat{A}$-genus of $TP/B$, and $\Ch_G \left(\ce/\cs \right)$ the
equivariant twisting Chern character of $\ce$, compare
\cite{bott78}. Hence on the level of cohomology the map
$\widetilde{\pi}^*$ defines a map $H^*_{\tau}(B_G) \to
H^*_{\tau}(P_G)$ which is pull-back composed with the multiplication
by $ \left(2 \pi i\right)^{-\frac{\dim P-\dim
B}{2}}\widehat{A}_{G}\left(TP/B\right) \Ch_{G} \left( \ce/ \cs
\right)$.

The second one is the map $\iota:\W_*(P^{(n)}) \to
CC^*(C_0^{\ify}(P^{(n)}))$ from the equation \eqref{iota}.

 The Proposition \ref{limit} immediately implies the
following:
\begin{cor}
The following  diagram commutes up to homotopy
\begin{equation}
\xymatrix{
&C^*\left(G, \W_*\left(P\right)\right)\ar[rr]^{\Phi\ }& & CC^*\left(C_0^{\ify}\left(P\right)\rtimes G\right)\\
&&&\\
&C^*\left(G,\W_*\left(B\right)\right)\ar[uu]^{\widetilde{\pi}^*}
\ar[rruu]^{\quad \Phi_{\Ab}}  & &
         }
\end{equation}

\end{cor}

Combining this result with the Theorem \ref{ltr} we obtain the main
result of this paper:

\begin{thm}
The following  diagram commutes up to homotopy
\begin{equation}
\xymatrix{
&C^*\left(G, \W_*\left(P\right)\right)\ar[rr]^{\Phi\ }& & CC^*\left(C_0^{\ify}\left(P\right)\rtimes G\right)\\
&&&\\
&C^*\left(G,\W_*\left(B\right)\right)\ar[uu]^{\widetilde{\pi}^*}
\ar[rr]^{\Phi}   & &CC^*\left(C_0^{\ify}\left(G\right)\right)
\ar[uu]^{\Ch\left(D\right)^t}
         }
\end{equation}
\end{thm}

\begin{rem}
In this paper we were concerned only with the smooth currents. Our
result can be generalized to the nonsmooth currents as follows.  Fix
arbitrary number $N$ and let $\W_k^N$ be the space of currents of
degree $k$ which belong to the Sobolev space $H^{N+k}_{\rm loc}$. We
can replace the complex $\W_*\left(B\right)$ by the complex using
these currents and all the results of these paper remain true. The
only change which needs to be made is in the estimates, as indicated
in the Remark \ref{differ}.
\end{rem}


\begin{thebibliography}{10}

\bibitem{bgv}
N.~Berline, E.~Getzler, and M.~Vergne.
\newblock {\em Heat kernels and {D}irac operators}, volume 298 of {\em
  Grundlehren der Mathematischen Wissenschaften [Fundamental Principles of
  Mathematical Sciences]}.
\newblock Springer-Verlag, Berlin, 1992.

\bibitem{bis85}
J.-M. Bismut.
\newblock The {A}tiyah-{S}inger index theorem for families of {D}irac
  operators: two heat equation proofs.
\newblock {\em Invent. Math.}, 83(1):91--151, 1985.

\bibitem{bf86}
J.-M. Bismut and D.~S. Freed.
\newblock The analysis of elliptic families. {II}. {D}irac operators, eta
  invariants, and the holonomy theorem.
\newblock {\em Comm. Math. Phys.}, 107(1):103--163, 1986.

\bibitem{bott78}
R.~Bott.
\newblock On some formulas for the characteristic classes of group-actions.
\newblock In {\em Differential topology, foliations and Gelfand-Fuks cohomology
  (Proc. Sympos., Pontif\'\i cia Univ. Cat\'olica, Rio de Janeiro, 1976)},
  volume 652 of {\em Lecture Notes in Math.}, pages 25--61. Springer, Berlin,
  1978.

\bibitem{bn}
J.-L. Brylinski and V.~Nistor.
\newblock Cyclic cohomology of \'etale groupoids.
\newblock {\em $K$-Theory}, 8(4):341--365, 1994.

\bibitem{cn83}
A.~Connes.
\newblock Cohomologie cyclique et foncteurs {${\rm Ext}\sp n$}.
\newblock {\em C. R. Acad. Sci. Paris S\'er. I Math.}, 296(23):953--958, 1983.

\bibitem{cn85}
A.~Connes.
\newblock Non-commutative differential geometry.
\newblock {\em Publ. Math. {IHES}}, 62:257--360, 1985.

\bibitem{cn86}
A.~Connes.
\newblock Cyclic cohomology and the transverse fundamental class of a
  foliation.
\newblock In {\em Geometric methods in operator algebras (Kyoto, 1983)}, volume
  123 of {\em Pitman Res. Notes Math. Ser.}, pages 52--144. Longman Sci. Tech.,
  Harlow, 1986.

\bibitem{cn94}
A.~Connes.
\newblock {\em Noncommutative Geometry}.
\newblock Academic Press, 1994.

\bibitem{cm90}
A.~Connes and H.~Moscovici.
\newblock Cyclic cohomology, the {N}ovikov conjecture and hyperbolic groups.
\newblock {\em Topology}, 29:345--388, 1990.

\bibitem{cm93}
A.~Connes and H.~Moscovici.
\newblock Transgression and the {C}hern character of finite-dimensional
  ${K}$-cycles.
\newblock {\em Comm. Math. Phys.}, 155(1):103--122, 1993.

\bibitem{cr99}
M.~Crainic.
\newblock Cyclic cohomology of \'etale groupoids: the general case.
\newblock {\em $K$-Theory}, 17(4):319--362, 1999.

\bibitem{crmo01}
M.~Crainic and I.~Moerdijk.
\newblock Foliation groupoids and their cyclic homology.
\newblock {\em Adv. Math.}, 157(2):177--197, 2001.

\bibitem{cuntz81}
J.~Cuntz.
\newblock Generalized homomorphisms between {$C\sp{\ast} $}-algebras and
  {$KK$}-theory.
\newblock In {\em Dynamics and processes (Bielefeld, 1981)}, volume 1031 of
  {\em Lecture Notes in Math.}, pages 31--45. Springer, Berlin, 1983.

\bibitem{cuntz82}
J.~Cuntz.
\newblock {$K$}-theory and {$C\sp{\ast} $}-algebras.
\newblock In {\em Algebraic $K$-theory, number theory, geometry and analysis
  (Bielefeld, 1982)}, volume 1046 of {\em Lecture Notes in Math.}, pages
  55--79. Springer, Berlin, 1984.

\bibitem{cuntz97}
J.~Cuntz.
\newblock Bivariante {$K$}-{T}heorie f\"ur lokalkonvexe {A}lgebren und der
  {C}hern-{C}onnes-{C}harakter.
\newblock {\em Doc. Math.}, 2:139--182 (electronic), 1997.

\bibitem{cq95a}
J.~Cuntz and D.~Quillen.
\newblock Cyclic homology and nonsingularity.
\newblock {\em J. Amer. Math. Soc.}, 8(2):373--442, 1995.

\bibitem{dup76}
J.~L. Dupont.
\newblock Simplicial de {R}ham cohomology and characteristic classes of flat
  bundles.
\newblock {\em Topology}, 15(3):233--245, 1976.

\bibitem{dup78}
J.~L. Dupont.
\newblock {\em Curvature and characteristic classes}.
\newblock Springer-Verlag, Berlin, 1978.
\newblock Lecture Notes in Mathematics, Vol. 640.

\bibitem{ft87}
B.~L. Fe{\u\i}gin and B.~L. Tsygan.
\newblock Cyclic homology of algebras with quadratic relations, universal
  enveloping algebras and group algebras.
\newblock In {\em $K$-theory, arithmetic and geometry (Moscow, 1984--1986)},
  pages 210--239. Springer, Berlin, 1987.

\bibitem{ge93a}
E.~Getzler.
\newblock Cartan homotopy formulas and the {G}auss-{M}anin connection in cyclic
  homology.
\newblock In {\em Quantum deformations of algebras and their representations
  (Ramat-Gan, 1991/1992; Rehovot, 1991/1992)}, volume~7 of {\em Israel Math.
  Conf. Proc.}, pages 65--78. Bar-Ilan Univ., Ramat Gan, 1993.

\bibitem{gj93}
E.~Getzler and J.~D.~S. Jones.
\newblock The cyclic homology of crossed product algebras.
\newblock {\em J. Reine Angew. Math.}, 445:161--174, 1993.

\bibitem{gs89}
E.~Getzler and Andr{\'a}s Szenes.
\newblock On the {C}hern character of a theta-summable {F}redholm module.
\newblock {\em J. Funct. Anal.}, 84(2):343--357, 1989.

\bibitem{gilkey95}
P.~B. Gilkey.
\newblock {\em Invariance theory, the heat equation, and the {A}tiyah-{S}inger
  index theorem}.
\newblock Studies in Advanced Mathematics. CRC Press, Boca Raton, FL, second
  edition, 1995.

\bibitem{goo85}
T.~G. Goodwillie.
\newblock Cyclic homology, derivations, and the free loopspace.
\newblock {\em Topology}, 24(2):187--215, 1985.

\bibitem{gor05}
A.~Gorokhovsky and J.~Lott.
\newblock {Local index theory over foliation groupoids}, arXiv:math.DG/0407246.

\bibitem{gor03}
A.~Gorokhovsky and J.~Lott.
\newblock Local index theory over \'etale groupoids.
\newblock {\em J. Reine Angew. Math.}, 560:151--198, 2003.

\bibitem{haefliger79}
A.~Haefliger.
\newblock Differential cohomology.
\newblock In {\em Differential topology (Varenna, 1976)}, pages 19--70.
  Liguori, Naples, 1979.

\bibitem{lg99b}
P.-Y. Le~Gall.
\newblock Th\'eorie de {K}asparov \'equivariante et groupo\"\i des. {I}.
\newblock {\em $K$-Theory}, 16(4):361--390, 1999.

\bibitem{lg99a}
P.-Y. Le~Gall.
\newblock Groupoid {$C\sp *$}-algebras and operator {$K$}-theory.
\newblock In {\em Groupoids in analysis, geometry, and physics (Boulder, CO,
  1999)}, volume 282 of {\em Contemp. Math.}, pages 137--145. Amer. Math. Soc.,
  Providence, RI, 2001.

\bibitem{lott92}
J.~Lott.
\newblock Superconnections and higher index theory.
\newblock {\em Geom. Funct. Anal.}, 2(4):421--454, 1992.

\bibitem{mw94a}
H.~Moscovici and F.~Wu.
\newblock Localization of topological {P}ontryagin classes via finite
  propagation speed.
\newblock {\em Geom. and Functional Analysis}, 1:52--92, 1994.

\bibitem{ni90}
V.~Nistor.
\newblock Group cohomology and the cyclic cohomology of crossed products.
\newblock {\em Invent. Math.}, 99(2):411--424, 1990.

\bibitem{ni91}
V.~Nistor.
\newblock A bivariant {C}hern character for $p$-summable quasihomomorphisms.
\newblock {\em $K$-Theory}, 5(3):193--211, 1991.

\bibitem{ni93a}
V.~Nistor.
\newblock A bivariant {C}hern-{C}onnes character.
\newblock {\em Ann. of Math. (2)}, 138(3):555--590, 1993.

\bibitem{ni95}
V.~Nistor.
\newblock Super-connections and non-commutative geometry.
\newblock In {\em Cyclic cohomology and noncommutative geometry (Waterloo, ON,
  1995)}, volume~17 of {\em Fields Inst. Commun.}, pages 115--136. Amer. Math.
  Soc., Providence, RI, 1997.

\bibitem{perd04}
D.~Perrot.
\newblock Retraction of the bivariant {C}hern character.
\newblock {\em $K$-Theory}, 31(3):233--287, 2004.

\bibitem{widom78}
H.~Widom.
\newblock Families of pseudodifferential operators.
\newblock In {\em Topics in functional analysis (essays dedicated to M. G.
  Kre\u\i n on the occasion of his 70th birthday)}, volume~3 of {\em Adv. in
  Math. Suppl. Stud.}, pages 345--395. Academic Press, New York, 1978.

\bibitem{widom80}
H.~Widom.
\newblock A complete symbolic calculus for pseudodifferential operators.
\newblock {\em Bull. Sci. Math. (2)}, 104(1):19--63, 1980.

\bibitem{wu97}
F.~Wu.
\newblock A bivariant {C}hern-{C}onnes character and the higher
  {$\Gamma$}-index theorem.
\newblock {\em $K$-Theory}, 11(1):35--82, 1997.

\end{thebibliography}
\end{document}